\documentclass[aap,MSNbibl,seceqn,citesort,dvips]{arximspdf}
\usepackage[vtexbibl]{}

\doi{10.1214/09-AAP656}
\volume{20}
\issue{4}
\pubyear{2010}
\firstpage{1470}
\lastpage{1511}

\makeatletter
\def\eqref#1{(\ref{#1})}
\newtheorem{maintheorem}{Theorem}
\newtheorem{theorem}{Theorem}[section]
\newtheorem{lemma}[theorem]{Lemma}
\newtheorem{claim}[theorem]{Claim}
\newtheorem{proposition}[theorem]{Proposition}
\newproclaim{observation}[theorem]{Observation}
\newtheorem{corollary}[theorem]{Corollary}
\newproclaim{remark}[theorem]{Remark}
\newproclaim{tttt}{Perfect allocation algorithm for $(1-\delta)n$ balls}
\newproclaim{bbb}{Basic version of allocation algorithm for $(1-\delta)n$ balls}
\newproclaim{iii}{Intermediate version of allocation algorithm for $(1-\delta)n$ balls}
\newcommand{\var}{\operatorname{Var}}
\makeatother

\begin{document}
\begin{frontmatter}

\title{Choice-memory tradeoff in allocations}
\runtitle{Choice-memory tradeoff in allocations}

\begin{aug}
\author[a]{\fnms{Noga} \snm{Alon}\thanksref{t1}\ead[label=e1]{nogaa@tau.ac.il}},
\author[b]{\fnms{Ori} \snm{Gurel-Gurevich}\ead[label=e2]{origurel@microsoft.com}\corref{}} \and
\author[b]{\fnms{Eyal} \snm{Lubetzky}\ead[label=e3]{eyal@microsoft.com}}
\runauthor{N. Alon, O. Gurel-Gurevich and E. Lubetzky}
\affiliation{Tel Aviv University, Microsoft Research and Microsoft Research}
\address[a]{N. Alon\\
School of Mathematics\\
Tel Aviv University\\
Tel Aviv, 69978\\ Israel\\
and\\
Microsoft-Israel R\&D Center\\
Herzeliya, 46725\\ Israel\\
\printead{e1}} 
\address[b]{O. Gurel-Gurevich\\E. Lubetzky\\
Microsoft Research\\
One Microsoft Way\\
Redmond, Washington 98052-6399\\ USA\\
\printead{e2}\\
\phantom{E-mail: }\printead*{e3}}
\thankstext{t1}{Supported in part by a USA Israeli BSF grant, by
a grant from
the Israel Science Foundation, by an ERC Advanced Grant
and by the Hermann Minkowski Minerva
Center for Geometry at Tel Aviv University.}
\end{aug}

\received{\smonth{10} \syear{2009}}

%
\begin{abstract}
In the classical balls-and-bins paradigm, where $n$ balls are placed
independently and uniformly in $n$ bins, typically
the number of bins with at least two balls in them is $\Theta(n)$ and
the maximum number of balls in a bin is $\Theta(\frac{\log n}{\log
\log n})$.
It is well known that when each round offers $k$ independent uniform
options for bins, it is possible to typically achieve a constant
maximal load
if and only if $k= \Omega(\log n)$. Moreover, it is possible w.h.p. to
avoid any collisions between
$n/2$ balls if $k> \log_2 n$.

In this work, we extend this into the setting where only $m$ bits of
memory are available.
We establish a tradeoff between the number of choices $k$ and the
memory $m$, dictated by the quantity $km/n$. Roughly put, we show that
for $k m \gg n$ one can achieve a constant maximal load, while for $k m
\ll n$ no substantial improvement can be gained over the case $k=1$
(i.e., a random allocation).

For any $k = \Omega(\log n)$ and $m=\Omega(\log^2 n)$, one can
achieve a constant load w.h.p. if $k m = \Omega(n)$, yet the
load is
unbounded if $km=o(n)$.
Similarly, if $k m > C n$ then $n/2$ balls can be allocated without any
collisions w.h.p.,
whereas for $k m < \varepsilon n$ there are typically $\Omega(n)$
collisions. Furthermore, we show that the load is w.h.p. at least
$\frac{\log(n/m)}{\log k + \log\log(n/m)}$.
In particular, for $k\leq\operatorname{polylog}(n)$, if $m =
n^{1-\delta}$ the
optimal maximal load is $\Theta(\frac{\log n}{\log\log n})$ (the
same as in the case $k=1$),
while $m=2n$ suffices to ensure a constant load.
Finally, we analyze nonadaptive allocation algorithms and give tight
upper and lower bounds for their performance.

\end{abstract}

%

\end{frontmatter}
%

\section{Introduction}\label{sec1}

The balls-and-bins paradigm (see, e.g., \cite{Feller,JK}) describes
the process where $b$ balls are placed independently and uniformly at
random in $n$ bins. Many variants of this classical occupancy problem
were intensively studied, having a wide range of applications in
computer science.

It is well known that when $b = \lambda n$ for $\lambda$ fixed and
$n\to\infty$, the load of each bin tends to Poisson with mean
$\lambda$ and the bins are asymptotically independent. In particular,
for $b=n$, the typical number of empty bins at the end of the process
is $(1/\mathrm{e}+o(1))n$. The typical maximal load in that case is
$(1+o(1))\frac{\log n}{\log\log n}$ (cf.~\cite{Gonnet}). In what
follows, we say that an event holds with high probability (w.h.p.) if
its probability tends to $1$ as $n\to\infty$.

The extensive study of this model in the context of load balancing was
pioneered by the celebrated paper of Azar et al. \cite{ABKU} (see
the survey \cite{MRS}) that analyzed the effect of a choice between
$k$ independent uniform bins on the maximal load, in an online
allocation of $n$ balls to $n$ bins. It was shown in \cite{ABKU} that
the \textsc{Greedy} algorithm (choose the least loaded bin of the
$k$) is
optimal and achieves a maximal-load of $\log_k\log n$ w.h.p., compared
to a load of $\frac{\log n}{\log\log n}$ for the original case $k=1$.
Thus, $k=2$ random choices already significantly reduce the maximal
load, and as $k$ further increases, the maximal load drops until it
becomes constant at \mbox{$k = \Omega(\log n)$}.

In the context of online bipartite matchings, the process of
dynamically matching each client in a group $A$ of size $n/2$ with one
of $k$ independent uniform resources in a group $B$ of size $n$
precisely corresponds to the above generalization of the balls-and-bins
paradigm: Each ball has $k$ options for a bin, and is assigned to one
of them by an online algorithm that should avoid collisions (no two
balls can share a bin). It is well known that the threshold for
achieving a perfect matching in this case is $k = \log_2 n$: For $k
\geq(1+\varepsilon)\log_2 n$, w.h.p. every client can
be exclusively matched to a target resource, and if $k \leq
(1-\varepsilon
)\log_2 n$ then $\Omega(n)$ requests cannot be satisfied.

In this work, we study the above models in the presence of a constraint
on the memory that the online algorithm has at its disposal. We find
that a tradeoff between the choice and the memory governs the ability
to achieve a perfect allocation as well as a constant maximal load.
Surprisingly, the threshold separating the subcritical regime from the
supercritical regime takes a simple form, in terms of the product of
the number of choices $k$, and the size of the memory in bits $m$:
\begin{itemize}
\item If $km \gg n$, then one can allocate $(1-\varepsilon)n$ balls in
$n$ bins without any collisions w.h.p., and consequently achieve
a load
of $2$ for $n$ balls.
\item If $km \ll n$, then \textit{any} algorithm for allocating
$\varepsilon n$ balls w.h.p. creates $\Omega(n)$ collisions
and an
unbounded maximal load.
\end{itemize}
Roughly put, when $km \gg n$ the amount of choice and memory at hand
suffices to guarantee an essentially best-possible performance. On the
other hand, when $km \ll n$, the memory is too limited to enable the
algorithm to make use
of the extra choice it has, and no substantial improvement can be
gained over the case $k=1$, where no choice is offered whatsoever.

Note that rigorous lower bounds for space, and in particular tradeoffs
between space and performance (time, communication, etc.), have been
studied intensively in the literature of Algorithm Analysis, and are
usually highly nontrivial. See, for example, \cite
{Ajtai,Beame,BJSK,BFMUW,BSSV,BC,Fortnow,FLMV} for some notable examples.

Our first main result establishes the exact threshold of the
choice-memory tradeoff for achieving a constant maximal-load. As
mentioned above, one can verify that when there is unlimited memory,
the maximal load is w.h.p. uniformly bounded iff $k = \Omega
(\log n)$.
Thus, assuming that $k=\Omega(\log n)$ is a prerequisite for
discussing the effect of limited memory on this threshold.

\begin{maintheorem}\label{thm-const-load}
Consider $n$ balls and $n$ bins, where each ball has $k=\Omega(\log
n)$ uniform choices for bins, and $m=\Omega(\log^2 n)$ bits of memory
are available.
If $k m = \Omega(n)$, one can achieve a maximal-load of $O(1)$ w.h.p.
Conversely, if $k m = o(n)$,
any algorithm w.h.p. creates a load that exceeds any constant.
\end{maintheorem}

Consider the case $k = \Theta(\log n)$. The na\"{\i}ve algorithm
for achieving a constant maximal-load in this setting requires roughly
$n$ bits of memory ($2n$ bits of memory always suffice; see
Section~\ref{subsec:m-size}). Surprisingly, the above theorem implies
that $O(n/\log n)$ bits of memory already suffice, and this is tight.

As we later show, one can extend the upper bound on the load, given in
Theorem~\ref{thm-const-load}, to $O(\frac{n}{km})$
(useful when $\frac{n}{km} \leq\frac{\log n}{\log\log n}$), whereas
the lower bound tends to~$\infty$ with $\frac{n}{km}$. This further
demonstrates how the quantity $\frac{n}{km}$ governs the value of the
optimal maximal load. Indeed, Theorem~\ref{thm-const-load} will follow
from Theorems \ref{thm-matching-lower} and~\ref{thm-matching-upper}
below, which
determine that the threshold for a perfect matching is $km = \Theta(n)$.

Again consider the case of $k = \Theta(\log n)$, where
an online algorithm with unlimited memory can achieve an $O(1)$ load
w.h.p. While the above theorem settles the memory threshold for
achieving a
constant load in this case, one can ask what the optimal maximal load
would be below the threshold. This is answered by the next theorem, which
shows that in this case, for example, $m=n^{1-\delta}$ bits of memory yield
no significant improvement over an algorithm which makes random allocations.

\begin{maintheorem}\label{thm-large-load}
Consider $n/k$ balls and $n$ bins, where each ball has $k$ uniform
choices for bins, and $m\geq\log n$ bits of memory are available.
Then for any algorithm, the maximal load is at least
$(1+o(1))\frac{\log(n/m)}{\log\log(n/m)+\log k}$ w.h.p.

In particular, if $m = n^{1-\delta}$ for some $\delta> 0$ fixed
and $2 \leq k \leq\operatorname{polylog}(n)$, then the maximal load
is $\Theta
(\frac{\log n}{\log\log n})$ w.h.p.
\end{maintheorem}

Recall that a load of order $\frac{\log n}{\log\log n}$ is what one
would obtain using a random allocation of $n$ balls in $n$ bins. The
above theorem states that, when $m=n^{1-\delta}$ and $k\leq
\operatorname{polylog}
(n)$, any algorithm would create such a load already after $n/k$ rounds.

Before describing our other results, we note that the lower bounds in
our theorems in fact apply to a more general setting. In the original
model, in each round the online algorithm chooses one of $k$ uniformly
chosen bins, thus inducing a distribution on the location of the next
ball. Clearly, this distribution has the property that no bin has a
probability larger than $k/n$.

Our theorems apply to a relaxation of the model, where the algorithm
is allowed to dynamically choose a distribution $Q_t$ for each round
$t$, which is required to satisfy the above property (i.e., $\|Q_t\|
_\infty\leq k/n$). We refer to these distributions as \textit{strategies}.

Observe that indeed this model gives more power to the online
algorithm. For instance, if $k=2$ (and the memory is unlimited), an
algorithm in the relaxed model can allocate $n/2$ balls perfectly (by
assigning $0$ probability to the occupied bins), whereas in the
original model collisions occur already with $n^{2/3}w(n)$ balls
w.h.p.,
for any $w(n)$ tending to $\infty$ with $n$.

Furthermore, we also relax the memory constraint on the model. Instead
of treating the algorithm as an automaton with $2^m$ states, we only impose
the restriction that there are at most $2^m$ different strategies to
choose from. In other words, at time~$t$, the algorithm knows the
entire history (the exact location of each ball so far), and needs to
choose one of its $2^m$ strategies for the next round. In this sense,
our lower bounds are for the case of limited communication complexity
rather than limited space complexity.

We note that all our bounds remain valid when each
round offers $k$ choices with repetitions.

\subsection{Tradeoff for perfect matching}

The next two theorems address the threshold for
achieving a perfect matching when allocating $(1-\delta) n$ balls in
$n$ bins for some fixed $0<\delta< 1$ [note that for $\delta= 0$,
even with unlimited memory, one needs $k = \Omega(n)$ choices to avoid
collisions w.h.p.].
The upper and lower bounds obtained for this threshold are tight up to
a multiplicative constant, and again
pinpoint its location at $k m = \Theta(n)$. The constants below were
chosen to simplify the proofs and could be optimized.

\begin{maintheorem}\label{thm-matching-lower}
For $\delta>0$ fixed, consider $(1-\delta) n $ balls and $n$ bins:
Each ball has $k$ uniform choices for bins, and there are $m\geq\log
n$ bits of memory.
If
\[
k m \leq\varepsilon n\qquad \mbox{for some small constant $\varepsilon>
0$} ,
\]
then \textit{any} algorithm has $\Omega(n)$ collisions w.h.p.

Furthermore, the maximal load is w.h.p. $\Omega(\log\log
(\frac{n}{km}))$.
\end{maintheorem}

\begin{maintheorem}\label{thm-matching-upper}
For $\delta>0$ fixed, consider $(1-\delta)n $ balls and $n$ bins,
where each ball has $k$ uniform choices for bins, and $m$ bits of
memory are available. The following holds for any $k \geq(3/\delta)
\log n$ and $m \geq\log n \cdot\log_2\log n$. If
\[
k m \geq C n\qquad \mbox{for some $C = C(\delta) > 0$},
\]
then a perfect allocation (no collisions) can be achieved w.h.p.
\end{maintheorem}

In light of the above, for any value of $k$, the online allocation
algorithm given by Theorem~\ref{thm-matching-upper} is optimal with
respect to its memory
requirements.

\subsection{Nonadaptive algorithms}

In the nonadaptive case, the algorithm is again allowed to choose a fixed
(possibly randomized) strategy for selecting the placement of
ball number $t$ in one of the $k$ possible randomly chosen bins
given in step~$t$. Therefore, each such algorithm consists of
a sequence $Q_1, Q_2, \ldots,Q_n$ of $n$ predetermined strategies,
where $Q_t$
is the strategy for selecting the bin in step number~$t$.

Here, we show
that even if $k= n\frac{\log\log n}{\log n}$, the maximum load
is w.h.p.  at least $(1-o(1)) \frac{\log n}{ \log
\log n}$,
that is, it is essentially as large as in the case $k=1$.
It is also possible\vspace*{1pt} to obtain tight bounds for larger values of $k$. We
illustrate this by considering the case $k=\Theta(n)$.

\begin{maintheorem}
\label{thm-nonadaptive}
Consider the problem of allocating $n$ balls into $n$ bins, where each
ball has $k$ uniform choices for bins,
using a nonadaptive algorithm.
\begin{longlist}[(ii)]
\item[(i)]\label{item-thm-nonadapt-1}
The maximum load in any nonadaptive algorithm with
$k \leq n\frac{\log\log n}{\log n}$ is w.h.p. at least $(1-o(1)) \frac{\log n}{ \log\log n}$.

\item[(ii)]\label{item-thm-nonadapt-2}
Fix $0 < \alpha< 1$. The maximum load in any nonadaptive algorithm with
$k=\alpha n$ is w.h.p.  $\Omega( \sqrt{\log
n})$. This is tight, that is, there exists a nonadaptive algorithm
with $k=\alpha n$ so that the maximum load in it is $O( \sqrt{\log
n})$ w.h.p.
\end{longlist}
\end{maintheorem}

\subsection{Range of parameters}\label{subsec:m-size}
In the above theorems and throughout the paper, the parameter $k$ may
assume values up to $n$. As for the memory, one may na\"{\i}vely use
$n\log_2 L$ bits
to store the status of $n$ bins, each containing at most $L$ balls.
The~next observation shows that the $\log_2 L$ factor is redundant.

\begin{observation*}\label{obs-memory-req}
At most $n+b-1$ bits of memory suffice to keep track of the number of
balls in each bin when allocating $b$ balls in $n$ bins.
\end{observation*}

Indeed, one can maintain the number of balls in each bin using a vector
in $\{0,1\}^{n+b-1}$, where $1$-bits stand for separators between the bins.
In light of this, the original case of unlimited memory corresponds to
the case $m = 2n$.

\subsection{Main techniques}
The key argument in the lower bound on the performance of the algorithm
with limited memory is analyzing the expected number of new collisions
that a given step introduces. We wish to estimate this value with an
error probability smaller than $2^{-m}$, so it would hold w.h.p. for
all of the $2^m$ possible strategies for this step.

To this end, we apply a large deviation inequality, which relates the
sum of a sequence of dependent random variables $(X_i)$ with the sum of
their ``predictions''~$(Y_i)$, where $Y_i$ is the expectation of $X_i$
given the history up to time $i$. Proposition~\ref{prop-predictions}
essentially shows that if the sum of the predictions $Y_i$ is large
(exceeds some $\ell$), then so is the sum of the actual random
variables $X_i$, except with probability $\exp(-c \ell)$. In the
application, the variable $X_i$ measures the number of new collisions
introduced by the $i$th ball, and $Y_i$ is determined by the strategy
$Q_i$ and the history so far.

The key ingredient in proving this proposition is a
Bernstein--Kolmogorov type inequality for martingales, which appears in
a paper of Freedman \cite{Freedman} from 1975, and bounds the
probability of deviation of a martingale in terms of
its cumulative variance. We reproduce its elegant proof for
completeness. Crucially, that theorem does not require a uniform bound on
individual variances (such as the one that appears in standard versions
of Azuma--Hoeffding), and rather treats them as random variables.
Consequently, the quality of our estimate in Proposition~\ref
{prop-predictions} is unaffected by the number of random variables involved.

For the upper bounds, the algorithm essentially partitions the bins
into blocks, where for different blocks it maintains an accounting
of the occupied bins with varying resolution. Once a block exceeds
a certain threshold of occupied bins, it is discarded and a new block
takes its place.

\subsection{Related work} The problem of balanced allocations with
limited memory is due to Itai Benjamini. In a recent independent work,
Benjamini and Makarychev \cite{BM} studied the special case of the
problem for $k=2$ (i.e., when there are two choices for bins at each
round). While our focus was mainly the regime
$k=\Omega(\log n)$ (where one can readily achieve a constant maximal
load when there is unlimited memory), our results
also apply for smaller values of $k$. Namely, as a by-product, we
improve the lower bound of \cite{BM} by a factor of $2$,
as well as extend it from $k=2$ to any $k\leq\operatorname{polylog}(n)$.

A different notion of memory was introduced to load balancing balls
into bins in \cite{MPS}, where one has the option of placing the
current ball in the least loaded bin offered in the previous round.
In that setting, one could indeed improve the asymptotics (yet not the
order) of the maximal load.
Note that in our case we consider the original balls and bins model (as
studied in \cite{ABKU}) and just impose restrictions on the space
complexity of the algorithm.

See, for example, \cite{MU}, Chapter 5, for more on the vast literature of
load balancing balls into bins and its applications in computer science.

A modern application for the classical online perfect matching problem
has advertisers (or \textit{bidders}) play the role of the bins and
internet search queries (or \textit{keywords}) play the role of the
balls. Upon receiving a search query, the search engine generates the
list of related advertisements (revealing the choices for this ball)
and must decide which of them to present in response to the query
(where to allocate the ball). Note that in the classical papers that
analyze online perfect matching one assumes a worst-case graph rather
than a random bipartite graph, and the requests are randomly permuted;
see \cite{KVV} for a fundamental paper in this area.

\subsection{Organization}
This paper is organized as follows. In
Section~\ref{sec:martingale}, we prove the large deviation inequality
(Proposition~\ref{prop-predictions}). Section~\ref{sec:lower-bound}
contains the lower bounds on the collisions and load, thus proving
Theorem~\ref{thm-matching-lower}. Section~\ref{sec:upper-bound}
provides algorithms for achieving a perfect-matching
and for achieving a constant load, respectively proving Theorem~\ref
{thm-matching-upper} and completing the proof of Theorem~\ref
{thm-const-load}. In Section~\ref{sec:large-q},
we extend the analysis of the lower bound to prove Theorem~\ref
{thm-large-load}. Section~\ref{sec:nonadaptive}
discusses nonadaptive allocations, and contains the proof of
Theorem~\ref{thm-nonadaptive}.
Finally, Section~\ref{sec:conclusion} is devoted to concluding remarks.

\section{A large deviation inequality}\label{sec:martingale}

This section contains a large deviation result,
which will later be one of the key ingredients in proving our lower
bounds for the load.
Our proof will rely on a Bernstein--Kolmogorov-type inequality of
Freedman \cite{Freedman}, which extends the standard Azuma--Hoeffding
martingale concentration inequality.
Given a sequence of bounded (possibly dependent) random variables
$(X_i)$ adapted to some filter $(\mathcal{F}_i)$,
one can consider the sequence $(Y_i)$ where $Y_i = \mathbb{E}
[X_i\mid
\mathcal{F}_{i-1}]$, which can be viewed as
predictions for the $(X_i)$'s. The following proposition essentially
says that, if the sum of the predictions $Y_i$ is large, so is the
sum of the actual variables $X_i$.

\begin{proposition}\label{prop-predictions}
Let $(X_i)$ be a sequence of random variables adapted to the filter
$(\mathcal{F}_i)$ so that $0 \leq X_i \leq M$ for all $i$, and let
$Y_i = \mathbb{E}
[X_i\mid\mathcal{F}_{i-1}]$. Then
\[
\mathbb{P}\biggl(\biggl\{\bigg|\frac{\sum_{i\leq t} X_i}{\sum
_{i\leq t}
Y_i}-1\bigg| \geq\frac12 \mbox{ and } \sum_{i\leq t} Y_i
\geq h\biggr\}\mbox{ for some $t$}\biggr)\leq
\exp\biggl(-\frac{h}{20M} + 2\biggr).
\]
\end{proposition}

\begin{pf}
As mentioned above, the proof hinges on a tail-inequality for sums of
random variables, which appears
in the work of Freedman \cite{Freedman} from 1975 (see also~\cite
{Steiger}), and extends such inequalities of Bernstein and Kolmogorov
to the setting of martingales. See \cite{Freedman} and the references
therein for more background on these inequalities, as well as \cite{Burkholder}
for similar martingale estimates.
We include the short proof of Theorem~\ref{thm-martingale} for completeness.

\begin{theorem}[(\cite{Freedman}, {Theorem 1.6})]\label{thm-martingale}
Let $(S_0,S_1,\ldots)$ be a martingale with respect to the filter
$(\mathcal{F}_i)$.
Suppose that $S_{i+1}-S_i \leq M$ for all $i$, and write $V_t = \sum
_{i=1}^t \var(S_i\mid\mathcal{F}_{i-1})$.
Then for any $s,v>0$ we have
\[
\mathbb{P}(S_n \geq S_0 + s , V_n \leq v\mbox{ for some
$n$})
\leq\exp\biggl[- \frac{s^2}{2(v + M s )}\biggr].
\]
\end{theorem}

\begin{pf}
Without loss of generality, suppose $S_0=0$, and put $X_i \stackrel
{\triangle}{=}S_i -
S_{i-1}$. Re-scaling $S_n$ by $M$, it clearly suffices to treat the
case $X_i \leq1$. Set
\[
V_t \stackrel{\triangle}{=}\sum_{i=1}^t \var(S_i\mid\mathcal
{F}_{i-1}) = \sum_{i=1}^t \mathbb{E}
(X_i^2 \mid\mathcal{F}_{i-1}),
\]
and for some $\lambda> 0$ to be specified later, define
\[
Z_t \stackrel{\triangle}{=}\exp\bigl(\lambda S_t - (\mathrm
{e}^\lambda- 1 - \lambda
)V_t\bigr).
\]
The next calculation will show that $(Z_t)$ is a super-martingale with
respect to the filter $(\mathcal{F}_t)$.
First, notice that the function
\[
f(z) \stackrel{\triangle}{=}\frac{\mathrm{e}^{z}-1-z}{z^2}\qquad \mbox{for $z \neq
0$}, f(0)\stackrel{\triangle}{=}\frac12,
\]
is monotone increasing [as $f'(z) > 0$ for all $z \neq0$], and in
particular, $ f(\lambda z) \leq f(\lambda)$
for all $ z \leq1$. Rearranging,
\[
\exp(\lambda z) \leq1 + \lambda z + (\mathrm{e}^{\lambda} - 1
- \lambda)z^2
\qquad \mbox{ for all $z \leq1$}.
\]
Now, since $X_i\leq1$ and $\mathbb{E}[X_i\mid\mathcal{F}_{i-1}]=0$
for all $i$, it
follows that
\begin{eqnarray*}
\mathbb{E}[ \exp(\lambda X_i) \mid\mathcal{F}_{i-1}]
&\leq&
1 + (\mathrm{e}^{\lambda} - 1 - \lambda)\mathbb
{E}[X_i^2
\mid\mathcal{F}_{i-1}] \\
&\leq&\exp\bigl((\mathrm{e}^{\lambda} - 1 - \lambda
)\mathbb{E}[X_i^2 \mid\mathcal{F}_{i-1}]\bigr).
\end{eqnarray*}
By definition, this precisely says that $\mathbb{E}[Z_i \mid\mathcal
{F}_{i-1}] \leq
Z_{i-1}$. That is, $(Z_t)$ is a super-martingale,
and hence by the Optional Stopping Theorem so is $(Z_{\tau\wedge n})$,
where $n$ is some integer and $\tau= \min\{t \dvtx  S_t \geq s\}$. In particular,
\[
\mathbb{E}Z_{\tau\wedge n} \leq Z_0 = 1,
\]
and (noticing that $V_{t+1} \geq V_t$ for all $t$) Markov's inequality
next implies that
\begin{eqnarray*}
\mathbb{P}\biggl(\bigcup_{t\leq n}(S_t \geq s,\ V_t \leq v)\biggr)
\leq\exp
[-\lambda s + (\mathrm{e}^{\lambda}-1-\lambda)v ].
\end{eqnarray*}
A choice of $\lambda= \log(\frac{s+v}v) \geq\frac
{s}{s+v} + \frac12 (\frac{s}{s+v})^2$ therefore yields
\begin{eqnarray*}
\mathbb{P}\biggl(\bigcup_{t\leq n}(S_t \geq s ,\ V_t \leq v)
\biggr)&\leq&
\exp\biggl[s-(s+v)\log\biggl(\frac{s+v}{v}\biggr)\biggr]\\
& \leq&\exp
\biggl[-\frac{s^2}{2(s+v)}\biggr],
\end{eqnarray*}
and taking a limit over $n$ concludes the proof.
\end{pf}

\begin{remark*}
Note that Theorem~\ref{thm-martingale} generalizes the
well-known version of the Azuma--Hoeffding inequality, where each of the
terms $\var(X_{i}\mid\mathcal{F}_{i-1})$ is bounded by some \textit{constant}
$\sigma_i^2$ (cf., e.g., \cite{McDiarmid}).
\end{remark*}


We now wish to infer Proposition~\ref{prop-predictions} from
Theorem~\ref{thm-martingale}. To this end, define
\[
Z_t \stackrel{\triangle}{=}\sum_{i=1}^t Y_i - X_i ,\qquad   V_t
\stackrel{\triangle}{=}\sum
_{i=1}^t\var( Z_i\mid\mathcal{F}_{i-1}),
\]
and observe that $(Z_t)$ is a martingale by the definition $Y_i =
\mathbb{E}
[X_i\mid\mathcal{F}_{i-1}]$. Moreover, as the $X_i$'s are uniformly bounded,
so are the increments of $Z_t$:
\[
| Z_i - Z_{i-1} | = | Y_i - X_i | \leq M.
\]
Furthermore, crucially, the variances of the increments are bounded as
well in terms
of the conditional expectations:
\[
\var(Y_i - X_i \mid\mathcal{F}_{i-1}) = \var(X_i
\mid
\mathcal{F}_{i-1}) \leq M \cdot\mathbb{E}[
X_i \mid\mathcal{F}_{i-1}] = M \cdot Y_i,
\]
giving that $V_t \leq M \sum_{i=1}^t Y_i$.

Finally, for any integer $j\geq1$ let $A_j$ denote the event
\[
A_j = \biggl( \biggl\{\sum_{i \leq t} X_i \leq\frac12 \sum
_{i\leq t}Y_i \mbox{ and }jh \leq\sum_{i \leq t} Y_i \leq(j+1)h
\biggr\}\mbox{ for some $t$}\biggr).
\]
%
Note that the event $A_j$ implies that $Z_t \geq jh/2$. Hence, applying
Theorem~\ref{thm-martingale} to the martingale $(Z_t)$ along with its
cumulative variances
$(V_t)$ we now get
\begin{eqnarray*}
\mathbb{P}(A_j) &\leq&\mathbb{P}\biggl( Z_t \geq\frac
{j}2 h  ,  V_t
\leq(j+1)hM\biggr)\\
& \leq&\exp\biggl[-\frac{(({j}/2) h)^2}{2
((j+1)hM+M(({j}/2)h))}\biggr]\\
&=& \exp\biggl[-\frac{j^2}{4(3j+2)}(h/M)\biggr] \leq\exp
\biggl(-\frac{h}{20M} j\biggr).
\end{eqnarray*}
Summing over the values of $j$, we obtain that if $h \geq20M$ then
\[
\mathbb{P}\biggl(\bigcup_{j\geq1} A_j\biggr) \leq\frac{\mathrm
{e}}{\mathrm
{e}-1} \exp\biggl(-\frac{h}{20M} \biggr) \leq\exp\biggl(-\frac
{h}{20M} j + 1\biggr),
\]
while for $h\leq20M$ the above inequality holds trivially. Hence, for
all $h>0$,
%
\begin{equation}\label{eq-prob-x/y-1}
\quad \mathbb{P}\biggl(\exists t \dvtx  \biggl\{\sum_{i \leq t} X_i \leq\frac
12 \sum_{i\leq t}Y_i \mbox{ and }\sum_{i \leq t} Y_i \geq h \biggr\}
\biggr)
\leq\exp\biggl(-\frac{h}{20M}+1\biggr).
\end{equation}
To complete the proof of the proposition, we repeat the above analysis for
\[
Z'_t \stackrel{\triangle}{=}-Z_t = \sum_{i=1}^t X_i - Y_i
,\qquad  V'_t \stackrel{\triangle}{=}\sum
_{i=1}^t\var( Z'_i\mid\mathcal{F}_{i-1}) = V_t.
\]
Clearly, we again have $|Z'_i - Z'_{i-1}| \leq M$ and $V'_t \leq M \sum
_{i=1}^t Y_i$. Defining
\[
A'_j = \biggl( \biggl\{\sum_{i \leq t} Y_i \leq\frac23 \sum
_{i\leq t}X_i \mbox{ and }jh \leq\sum_{i \leq t} Y_i \leq(j+1)h
\biggr\}\mbox{ for some $t$}\biggr),
\]
it follows that the event $A'_j$ implies that $Z'_t \geq\frac12 \sum
_{i}Y_i \geq jh/2$. Therefore, as before, we have that
\[
\mathbb{P}(A'_j) \leq\mathbb{P}\biggl( Z'_t \geq\frac
{j}2h  ,  V'_t
\leq(j+1)hM\biggr) \leq\exp\biggl(-\frac{h}{20M} j\biggr),
\]
and thus for all $h>0$
%
\begin{equation}\label{eq-prob-x/y-2}
\mathbb{P}\biggl(\exists t\dvtx  \biggl\{\sum_{i \leq t} Y_i \leq\frac
23 \sum_{i\leq t}X_i \mbox{ and }\sum_{i \leq t} Y_i \geq h \biggr\}
\biggr)
\leq\exp\biggl(-\frac{h}{20M}+1\biggr).
\end{equation}
Summing the probabilities in \eqref{eq-prob-x/y-1} and \eqref
{eq-prob-x/y-2} yields the desired result.
\end{pf}

We note that essentially the same proof yields the following
generalization of Proposition~\ref{prop-predictions}.
As before, the constants can be optimized.

\begin{proposition}\label{prop-predictions-2}
Let $(X_i)$ and $(Y_i)$ be as given in Proposition~\textup{\ref
{prop-predictions}}. Then for any $0 < \varepsilon\leq\frac12$,
\begin{eqnarray*}
&&\mathbb{P}\biggl(\biggl\{\bigg|\frac{\sum_{i\leq t} X_i}{\sum
_{i\leq t}
Y_i} -1 \bigg| \geq\varepsilon\mbox{ and } \sum_{i\leq t} Y_i \geq
h\biggr\}\mbox{ for some $t$}\biggr)\\
&&\qquad  \leq\exp\biggl(-\frac
{h\varepsilon
^2}{5M} + 2\biggr).
\end{eqnarray*}
%
\end{proposition}

\begin{remark}\label{rem-prop-stopping-time}
The statements of Propositions~\ref{prop-predictions} and~\ref{prop-predictions-2} hold also in conjunction with any
stopping time $\tau$ adapted to the filter $(\mathcal{F}_i)$. That
is, we get
the same bound on the probability of the mentioned event happening at
any time $t < \tau$. This follows easily, for instance, by altering
the sequence of increments to be identically $0$ after $\tau$.
Such statements become useful when the uniform bound on the increments
is only valid before $\tau$.
\end{remark}

\section{Lower bounds on the collisions and load}\label{sec:lower-bound}

In this section, we prove Theorem~\ref{thm-matching-lower} as well as
the lower bound in Theorem~\ref{thm-const-load}, by showing that if
the quantity $km/n$
is suitably small, then any allocation would necessarily produce nearly
linearly many bins with arbitrarily large load.

The main ingredient in the proof is a bound for the number of
collisions, that is, pairs of balls that share a bin, defined next.
Let $N_t(i)$ denote the number of balls in bin $i$ after performing $t$
rounds; the number of collisions at time $t$ is then
\[
\operatorname{Col}_2(t) \stackrel{\triangle}{=}\sum_{i=1}^n \pmatrix{ N_t(i)\cr 2}.
\]
The following theorem provides a lower bound on $\operatorname
{Col}_2(t)$ for $t
\geq c\cdot km$
for some absolute $c > 0$.

\begin{theorem}\label{thm-collisions-lower-bound}
Consider $n $ balls and $n$ bins,
where each ball has $k$ uniform choices for bins, and $m\geq\log n$
bits of memory are available.
%
\begin{longlist}[(ii)]
\item[(i)]\label{item-thm-col-1}For all $t \geq500\cdot k m $, we have
\[
\mathbb{E}\operatorname{Col}_2(t) \geq t^2 / (9n).
\]
\item[(ii)]\label{item-thm-col-2}Furthermore, with probability
$1-O(n^{-4})$, for all $L=L(n)$ and any $t \geq(  500 k m
\vee  30 \sqrt{Ln\log n} )$, either the maximal load is at
least $L$ or
\[
\operatorname{Col}_2(t) \geq t^2 / (16n).
\]
\end{longlist}
\end{theorem}

Note that the main statement of Theorem~\ref{thm-matching-lower}
immediately follows from the above theorem,
by choosing $t = (1-\delta)n$ and $L = \sqrt{n}$. Indeed, recalling
the assumption in Theorem~\ref{thm-matching-lower} that $m \geq\log
n$, we obtain
that, except with probability $O(n^{-4})$,
either the algorithm creates a load of $\sqrt{n}$, or
it has $\operatorname{Col}_2(n) \geq\frac{(1-\delta)^2}{16} n$.
Observing that a
load of $L$ immediately induces ${L\choose 2}$ collisions, we deduce
that either way there are at least $\Omega(n)$ collisions w.h.p.

We next prove Theorem~\ref{thm-collisions-lower-bound}; the
statement of Theorem~\ref{thm-matching-lower} on unbounded maximal
load will follow from an iterative
application of a more general form of this theorem (namely,
Theorem~\ref{thm-collisions-gen}), which appears
in Section~\ref{subsec:unbounded-maxload}.

\begin{pf*}{Proof of Theorem~\ref{thm-collisions-lower-bound}}
As noted in the \hyperref[sec1]{Introduction}, we relax the model by allowing the
algorithm to choose any distribution $\mu= (\mu(1),\ldots,\mu(n))$
for the location of the next ball, as long as it satisfies
$\|\mu\|_\infty\leq k/n$.

We also relax the memory constraint as follows. The algorithm has a
pool of at most $2^m$ different strategies, and may choose any of them
at a given step without any restriction (basing its dynamic decision on
the entire history).

To summarize, the algorithm has a pool of at most $2^m$ strategies, all
of which have an $L^\infty$-norm of at most $k/n$. In each given
round, it adaptively chooses a strategy $\mu$ from this pool based on
the entire history, and a ball then falls to a bin distributed
according to $\mu$.

The outline of the proof is as follows: consider the sequence
$Q_1,\ldots,Q_n$, chosen adaptively out of the pool of $2^m$ of strategies.
The large deviation inequality of Section~\ref{sec:martingale}
(Proposition~\ref{prop-predictions})
will enable us to show the following. The expected number of collisions
encountered in the above process
is well approximated by the expected number of collisions between $n$
independent balls, placed according to $Q_1,\ldots,Q_n$
(i.e., equivalent to the result of the nonadaptive algorithm with
strategies $Q_1,\ldots,Q_n$).

Having reduced the problem to the analysis of a nonadaptive algorithm,
we may then
derive a lower bound on $\mathbb{E}\operatorname{Col}_2(t)$ by
analyzing the structure of
the above strategies.
This bound is then translated to a bound on $\operatorname{Col}_2(t)$
using another
application
of the large deviation inequality of Proposition~\ref{prop-predictions}.

Let $\nu= (\nu(1),\ldots,\nu(n))$ be an arbitrary probability
distribution on $[n]$ satisfying
$\|\nu\|_\infty\leq k/n$, and denote by $Q_s = (Q_s(1),\ldots
,Q_s(n))$ the
strategy of the algorithm at time $s$. It will be convenient from time
to time to treat these distributions as vectors in $\mathbb{R}^n$.

By the above discussion, $Q_s$ is a random variable
whose values belong to some a priori set $\{\mu_1,\ldots
,\mu_{2^m}\}$.
We further let $J_s$ denote the actual position of the ball at time $s$
(drawn according
to the distribution $Q_s$).

Given the strategy at time $s$, let $x_s$ denote the probability of a
collision between
$\nu$ and $Q_s$ given $J_s$, that is, that the ball that is distributed
according to $\nu$
will collide with the one that arrived in time $s$. We let $v_s$ be the
inner product of $Q_s$ and~$\nu$, which
measures the expectation of these collisions:
\begin{eqnarray*}
x_s^\nu&\stackrel{\triangle}{=}&\nu( J_s),\\
v_s^{\nu} &\stackrel{\triangle}{=}&\langle Q_s,\nu\rangle  = \sum
_{i=1}^n Q_s(i) \nu(i)
= \mathbb{E}[x_s^\nu\mid\mathcal{F}_{s-1}].
\end{eqnarray*}
Further define the cumulative sums of $v_s^\nu$ and $x_s^\nu$ as follows:
\begin{eqnarray*}
X_t^{\nu} &\stackrel{\triangle}{=}&\sum_{s=1}^t x_s^\nu,\\
V_t^{\nu} &\stackrel{\triangle}{=}&\sum_{s=1}^t v_s^\nu.
\end{eqnarray*}
To motivate these definitions, notice that given the history up to time
$s-1$ and any possible strategy for the next round, $\nu$,
we have
\[
X_{s-1}^{\nu} = \sum_{i=1}^{s-1} \nu(J_i) = \sum_{i=1}^n \nu(i) |\{
r<s\dvtx J_r = i\}|
= \sum_{i=1}^n \nu(i) N_{s-1}(i),
\]
and so $X_{s-1}^{Q_s}$ is the expected number of collisions that will
be contributed by the ball $J_s \sim Q_s$
given the entire history $\mathcal{F}_{s-1}$. Summing over $s$, we
have that
\[
\mathbb{E}\operatorname{Col}_2(t) = \mathbb{E}\Biggl[\sum_{s=1}^t
X_{s-1}^{Q_s}\Biggr],
\]
thus estimating the quantities $X_{s-1}^{Q_s}$ will provide a bound on
the expected number of collisions.
Our aim in the next lemma is to show that w.h.p., whenever
$V_{s-1}^{Q_s}$ is large, so is~$X_{s-1}^{Q_s}$. This will reduce
the problem to the analysis of the quantities~$V_{s-1}^{Q_s}$, which
are deterministic functions of~$Q_1,\ldots,Q_n$.
This is the main conceptual ingredient in the lower bound, and its
proof will follow directly
from the large deviation estimate given in Proposition~\ref{prop-predictions}.

\begin{lemma}\label{lem-Xt-Wt-approx}
Let $Q_1,\ldots,Q_n$ be a sequence of strategies adapted to the filter~$(\mathcal{F}_{i})$, and let $X^\nu_s$ and $V^\nu_s$ be defined as above.
Then with probability at least $1-O(\mathrm{e}^{-4m})$, for every $\nu
\in\{\mu_1,\ldots,\mu_{2^m}\}$ and every $s$
we have that $V^\nu_s \geq100\|\nu\|_\infty m$ implies $X^\nu_s
\geq V^\nu_s/2$.
\end{lemma}

\begin{pf}
Before describing the proof, we wish to emphasize a delicate point.
The lemma holds for any sequence of strategies $Q_1,Q_2,\ldots,Q_n$
(each $Q_i$ is an arbitrary function of $\mathcal{F}_{i-1}$). No restrictions
are made here on the way each
such $Q_i$ is produced (e.g., it does not even need to belong to the
pool of $2^m$ strategies), as long as it satisfies $\|Q_i\|_\infty\leq k/n$.
The reason that such a general statement is possible is the following:
Once we specify how each $Q_i$ is determined from $\mathcal{F}_{i-1}$ (this
can involve extra random bits, in case the adaptive algorithm is
randomized), the process of exposing the positions
of the balls, $J_i \sim Q_i$, defines a martingale. Hence, for each
fixed $\nu$, we would be able to show that the desired event occurs except
with probability $O(\mathrm{e}^{-5m})$. A union bound over the
strategies $\nu$ (which, crucially, do belong to the pool of size $2^m$)
will then complete the proof.

Fix a strategy $\nu$ out of the pool of $2^m$ possible strategies, and
recall the definitions of $x_s^\nu$ and $v_s^\nu$, according to which
\[
0 \leq x_s^\nu\leq\|\nu\|_\infty,\qquad   v_s^{\nu} = \mathbb
{E}[x_s^\nu
\mid\mathcal{F}_{s-1}].
\]
By applying Proposition~\ref{prop-predictions} to the sequence
$(x_s^\nu)$ (with the cumulative sums $X_s^\nu$ and cumulative
conditional expectations $V_s^\nu$), we obtain that for all $h$,
\[
\mathbb{P}( X_s^\nu\leq V_s^\nu/2, V_s^\nu\geq h\mbox{ for some $s$})
\leq\exp\biggl(-\frac{h}{20\|\nu\|_\infty} + 2\biggr).
\]
Thus, taking $h = 100\|\nu\|_\infty m$ we obtain that
\[
\mathbb{P}( X_s^\nu\leq V_s^\nu/2, V_s^\nu\geq100\|\nu\|_\infty
m\mbox
{ for some $s$}) \leq\exp(-5m+2).
\]
Summing over the pool of at most $2^m$ predetermined strategies, $\nu$
completes the proof.
\end{pf}

Having shown that $X_t^\nu$ is well approximated by $V_t^\nu$, and
recalling that we are interested
in estimating $X_{s-1}^{Q_s}$, we now turn our attention to the
possible values of~$V_{s-1}^{Q_s}$.

\begin{claim}\label{clm-large-Ws-sum}
For any sequence of strategies $Q_1,\ldots,Q_t$, we have that
\[
\sum_{s=1}^{t} V_{s-1}^{Q_s} \geq\frac{t(t-k)}{2n}.
\]
\end{claim}

\begin{pf}
By our definitions, for the strategies $Q_1,\ldots,Q_t$ we have
\begin{eqnarray}\label{eq-Ws-sum-estimate}
\sum_{s=1}^t V_{s-1}^{Q_s} &=& \sum_{s=1}^t \sum_{r=1}^{s-1}
\langle Q_r,Q_s\rangle  = \sum_{i=1}^n \sum_{r < s \leq t} Q_r(i) Q_s(i)
\nonumber\\[-8pt]\\[-8pt]
&=& \frac12 \sum_{i=1}^n \Biggl[\Biggl(\sum_{s=1}^t Q_s(i)\Biggr)^2 -
\sum_{s=1}^t Q_s(i)^2\Biggr].\nonumber
\end{eqnarray}
Recalling the definition of the strategies $Q_i$, we have that
\[
\cases{
0 \leq Q_s(i) \leq k/n, & \quad for all $i$ and $s$,\vspace*{2pt}\cr
\displaystyle\sum_{i=1}^n Q_s(i) = 1, &\quad for all $s$.
}
\]
Therefore,
\[
\sum_{i=1}^n\sum_{s=1}^t Q_s(i)^2 \leq\frac{k}n \sum_{i=1}^n \sum
_{s=1}^t Q_s(i) = \frac{k t}n.
\]
On the other hand, by Cauchy--Schwarz,
\[
\sum_{i=1}^n\Biggl(\sum_{s=1}^t Q_s(i)\Biggr)^2 \geq\frac{1}n
\Biggl(\sum_{i=1}^n \sum_{s=1}^t Q_s(i)\Biggr)^2 = \frac{t^2}n.
\]
Plugging these two estimates in \eqref{eq-Ws-sum-estimate}, we deduce that
\[
\sum_{s=1}^t V_{s-1}^{Q_s} \geq\frac{t(t-k)}{2n},
\]
as required.
\end{pf}

While the above claim tells us that the average size of $V_{s-1}^{Q_s}$
is fairly large [has order at least $(t-k)/n$], we wish to obtain
bounds corresponding to individual distributions $Q_s$.
As we next show, this sum indeed enjoys a significant contribution from
indices $s$ where $V_{s-1}^{Q_s} = \Omega(km/n)$.
More precisely, setting $h = 100km/n$, we claim that
for large enough $n$,
%
\begin{equation}\label{eq-many-large-Ws}
\sum_{s=1}^t V_{s-1}^{Q_s} \bolds{1}_{\{V_{s-1}^{Q_s} > h\}} \geq
\frac
{t^2}{4n} .
\end{equation}
To see this, observe that if
\[
t \geq t_0 \stackrel{\triangle}{=}5hn = 500 k m,
\]
then
\[
\sum_{s=1}^t V_{s-1}^{Q_s} \bolds{1}_{\{V_{s-1}^{Q_s} \leq h\}} \leq
t h
\leq\frac{t^2}{5n}.
\]
Combining this with Claim~\ref{clm-large-Ws-sum} [while noting that
$\frac{t(t-k)}{2n} = (1-o(1))\frac{t^2}{2n}$] yields \eqref
{eq-many-large-Ws} for any sufficiently large $n$.

We may now apply Lemma~\ref{lem-Xt-Wt-approx}, and obtain that, except
with probability $O(\mathrm{e}^{-4m})$, whenever
$ V_{s-1}^{Q_s} > h$ we have $X_{s-1}^{Q_s} \geq\frac12
V_{s-1}^{Q_s}$, and so
%
\begin{equation}\label{eq-Xs-sum-t^2/8n}
\sum_{s=1}^t X_{s-1}^{Q_s}\geq\frac12 \sum_{s=1}^t V_{s-1}^{Q_s}
\bolds{1}_{\{V_{s-1}^{Q_s} > h\}} \geq\frac{t^2}{8n} \qquad \mbox{for
all $t
\geq t_0$}.
\end{equation}
Altogether, since $X_{s-1}^{Q_s} \geq0$, we infer that
%
\begin{equation}
\label{eq-E-Col(t)}
\quad \mathbb{E}\operatorname{Col}_2(t) = \mathbb{E}\Biggl[\sum_{s=1}^t
X_{s-1}^{Q_s}\Biggr] \geq\frac
{t^2}{8n}\bigl(1- O(n^{-4})\bigr) \geq\frac{t^2}{9n}\qquad \mbox{for all
$t \geq t_0$},
\end{equation}
where the last inequality holds for large enough $n$. This proves
part~(i) of Theorem~\ref{thm-collisions-lower-bound}.

It remains to establish concentration for $\operatorname{Col}_2(t)$
under the
additional assumption that $t \geq30 \sqrt{L n\log n}$
for some $L = L(n)$. First, set the following stopping-time for
reaching a maximal-load of $L$:
\[
\tau_L \stackrel{\triangle}{=}\min\Bigl\{t \dvtx  \max_j N_t(j) \geq
L\Bigr\}.
\]
Next, recall that
\[
\operatorname{Col}_2(t) = \sum_{s=1}^t N_{s-1}(J_{s}),
\]
and notice that
\[
 \mathbb{E}[ N_{s-1}(J_s) \mid\mathcal{F}_{s-1}] =
\sum_{i=1}^n
Q_s(i)N_{s-1}(i) = X_{s-1}^{Q_s}.
\]
Therefore, we may apply our large deviation estimate given in
Section~\ref{sec:martingale} (Proposition~\ref{prop-predictions}),
combined with the stopping-time $\tau_L$ (see Remark~\ref
{rem-prop-stopping-time}):
\begin{itemize}
\item The sequence of increments is $(N_{s-1}(J_s))$.
\item The sequence of conditional expectations is $(X_{s-1}^{Q_s})$.
\item The bound on the increments is $L$, as $N_{s-1}(J_s) \leq\max_i
N_{s-1}(i) \leq L$ for all $s < \tau_L$.
\end{itemize}
It follows that
\begin{eqnarray*}
&&\mathbb{P}\biggl(\biggl\{\operatorname{Col}_2(t) \leq\frac
12\sum_{s\leq t}
X_{s-1}^{Q_s} \mbox{ and } \sum_{s\leq t} X_{s-1}^{Q_s} \geq\frac
{t^2}{8n}\biggr\}\mbox{ for some $t < \tau_L$}\biggr) \\
&&\qquad \leq
\exp\biggl(-\frac{t^2/8n}{20L} + 2\biggr) \leq O(n^{-5}),
\end{eqnarray*}
where the last inequality is by the assumption $t \geq30\sqrt{Ln\log n}$.
Finally, by~\eqref{eq-Xs-sum-t^2/8n}, we also have
that $\sum_{s\leq t} X_{s-1}^{Q_s} \geq t^2/(8n)$ for all $t \geq
t_0$, except with probability $O(n^{-4})$.
Combining these two statements, we deduce that for any $t \geq(t_0
\vee 30\sqrt{L n\log n})$,
\[
\mathbb{P}\biggl(\operatorname{Col}_2(t) < \frac{t^2}{16n},~\tau
_L>t \biggr) = O
(n^{-4}),
\]
concluding the proof of Theorem~\ref{thm-collisions-lower-bound}.
\end{pf*}

\subsection{Boosting the subcritical regime to unbounded maximal
load}\label{subsec:unbounded-maxload}
While Theorem~\ref{thm-collisions-lower-bound} given above provides a
careful analysis for the number of $2$-collisions, that is, pairs of balls
sharing a bin, one can iteratively apply this theorem, with very few
modifications, in order to obtain that the number of $q$-collisions (a
set of $q$ balls sharing a bin) has order $\Omega(n^{1-o(1)})$ w.h.p.
The proof of this result hinges on Theorem~\ref{thm-collisions-gen}
below, which is a generalization of Theorem~\ref{thm-collisions-lower-bound}.

Recall that in the relaxed model studied so far, at any given time $t$
the algorithm adaptively selects a strategy $Q_t$ (based on the entire
history $\mathcal{F}_{t-1}$), after which a ball is positioned in a
bin $J_t
\sim Q_t$. We now introduce an extra set of random variables,
in the form of a sequence of increasing subsets, $A_1 \subset\cdots
\subset A_n \subset[n]$. The set $A_t$ is determined
by $\mathcal{F}_{t-1}$, and has the following effect: If $J_t \in
A_t$, we add
a ball to this bin as usual, whereas if
$J_t \notin A_t$, we ignore this ball (all bins remain unchanged). That
is, the number of balls in bin $i$ at time $t$ is now given by
\[
N_t(i) \stackrel{\triangle}{=}\sum_{s=1}^t \bolds{1}_{\{J_s = i\}}
\bolds{1}_{A_s}(i),
\]
and as before we are interested in a lower bound for the number of collisions:
\[
\operatorname{Col}_2(t) \stackrel{\triangle}{=}\sum_{i=1}^n \pmatrix{ N_t(i)\cr 2}.
\]
The idea here is that, in the application, the set $A_t$ will consists of
the bins that already contain $\ell$ balls at time $t$. As such, they indeed
form an increasing sequence of subsets determined by $(\mathcal
{F}_i)$. In
this case, any collision
corresponds to $2$ balls placed in some bin which already has $\ell$
other balls, and thus immediately implies a load of $\ell+2$.

\begin{theorem}\label{thm-collisions-gen}
Consider the following balls and bins setting:
\begin{enumerate}[(3)]
%
\item[(1)] The online adaptive algorithm has a pool of $2^m$ possible
strategies, where each strategy $\mu$ satisfies $\|\mu\|_\infty\leq k/n$.
The algorithm selects a (random) sequence of strategies $Q_1,\ldots
,Q_n$ adapted to the filter $(\mathcal{F}_i)$.
\item[(2)] Let $A_1\subset\cdots\subset A_n \subset[n] $ denote a random
increasing sequence of subsets adapted to the filter $(\mathcal
{F}_i)$, that is, $A_i$
is determined by $\mathcal{F}_{i-1}$.
\item[(3)] There are $n$ rounds, where in round $t$ a new potential location
for a ball is chosen according to $Q_t$. If this location
belongs to $A_t$, a ball is positioned there (otherwise, nothing happens).
\end{enumerate}
Define $T = \sum_{s=1}^n Q_s(A_s)$. Then for any $L=L(n)$,
\[
\mathbb{P}\biggl(T \geq30\bigl( \sqrt{k m n}  \vee\sqrt{L n \log
n}
\bigr),  \operatorname{Col}_2(n) < \frac{T^2}{16n}, \max_{j}N_n(j) \leq
L\biggr) \leq
O(n^{-4}).
\]
\end{theorem}

\begin{pf}
As the proof follows the same arguments of Theorem~\ref
{thm-collisions-lower-bound}, we restrict our attention to describing
the modifications that are required for the new statement to hold.

Define the following subdistribution of $Q_s$ with respect to $A_s$:
\[
Q'_s \stackrel{\triangle}{=}Q_s \bolds{1}_{A_s}.
\]
As before, given $Q_s$, the strategy at time $s$, define the following
parameters:
\[
x_s^\nu\stackrel{\triangle}{=}\nu( J_s),\qquad v_s^{\nu} \stackrel
{\triangle}{=}\sum_{i=1}^n Q'_s(i) \nu(i),
\]
and let the cumulative sums of $v_s^\nu$ and $x_s^\nu$ be denoted by
\[
X_t^{\nu} \stackrel{\triangle}{=}\sum_{s=1}^t x_s^\nu,\qquad V_t^{\nu
} \stackrel{\triangle}{=}\sum_{s=1}^t
v_s^\nu.
\]
We claim that a statement analogous to that of Lemma~\ref
{lem-Xt-Wt-approx} holds as is with respect to the above definitions,
for any choice of increasing subsets $A_1\subset\cdots\subset A_n$
[adapted to the filter $(\mathcal{F}_i)$]. As we soon argue, the martingale
concentration argument is valid without any changes, and the only
delicate point is the identity of the target strategy $\nu$.

\begin{lemma}\label{lem-Xt-Wt-approx-gen}
Let $Q_1,\ldots,Q_n$ and $A_1\subset\cdots\subset A_n$ be strategies and
subsets respectively, adapted to the filter $(\mathcal{F}_{i})$, and let
$X^\nu_s$
and $V^\nu_s$ be defined as above.
Then with probability at least $1-O(\mathrm{e}^{-4m})$, for every $\nu
\in\{\mu_1,\ldots,\mu_{2^m}\}$ and every~$s$
we have that $V^{\nu'}_s \geq100\|\nu\|_\infty m$ implies $X^{\nu
'}_s \geq V^{\nu'}_s/2$, where $\nu' = \nu\bolds{1}_{A_{s+1}}$.
\end{lemma}

\begin{pf}
Let $\nu$ be a strategy. Previously (in the proof of Lemma~\ref
{lem-Xt-Wt-approx}), we compared $X_s^\nu$ to $V_s^\nu$
using the large deviation inequality of Section~\ref{sec:martingale}.
Now, for each~$s$, our designated $\nu'$ is a function of $\nu$ and
$A_{s+1}$, and hence depends on $\mathcal{F}_s$.
In particular, there are potentially more than $2^m$ different
strategies to consider as $\nu'$,
destroying our union bound.
%
The crucial observation that resolves this issue is the following.

\begin{observation}
Let $r > s$ and let $\nu$ be a strategy. Then $V_{s}^\nu= V_{s}^{\nu
'}$ and $X_{s}^\nu= X_{s}^{\nu'}$
for any increasing sequence $A_1,\ldots,A_r$, where $\nu'=\nu\bolds
{1}_{A_r}$.
\end{observation}

To see this, first consider $X_s^\nu$ and $X_s^{\nu'}$. If $x_i^\nu$
for some $1 \leq i \leq s$ had a nonzero
contribution to $X_s^\nu$, then by definition \mbox{$J_i \in A_i$}.
Since~$A_i\subset A_r$, we also have \mbox{$J_i \in A_r$}, and so $x_i^{\nu'}
= \nu(J_i)\bolds{1}_{A_r}(J_i) = x_i^{\nu}$.
The statement $V_s^{\nu'}=V_s^{\nu}$ now follows from the fact that
$V_s^\nu$ is the sum of $v_i^\nu= \mathbb{E}[x_i^\nu\mid\mathcal
{F}_{i-1}]$.

Using the above observation, it now suffices to prove the statement of
Lemma~\ref{lem-Xt-Wt-approx-gen} directly on the strategies $\nu$
(rather than on $\nu'$). Hence, the only difference between this
setting and that of Lemma~\ref{lem-Xt-Wt-approx} is that
here some of the rounds are forfeited (as reflected in the new
definition of the $v_s^\nu$'s). The proof of
Lemma~\ref{lem-Xt-Wt-approx} therefore holds unchanged for this case.
\end{pf}

Similarly, the following claim is the analogue of Claim~\ref
{clm-large-Ws-sum}, with $t$ (the number of balls
in the original version) replaced by $T=\sum_s \sum_i Q'_s(i)$ (the
expected number of balls actually positioned).

\begin{claim}\label{clm-large-Ws-sum-gen}
For any $Q_1,\ldots,Q_n$ and $A_1\subset\cdots\subset A_n$, we have that
\[
\sum_{s=1}^{n} V_{s-1}^{Q'_s} \geq\frac{T(T-k)}{2n}.
\]
\end{claim}

The proof of the above claim follows from the exact same argument as in
Claim~\ref{clm-large-Ws-sum}.
Notice that the bound there, given as a function of $t$, was
actually a bound in terms of $\sum_{s=1}^t \sum_i Q_s(i)$, and so
replacing $Q_s$ by $Q'_s$ yields the desired bound as a function of $T$.

With this in mind, set $h = 100 km/n$ and note that, clearly,
\[
\sum_{s=1}^n V_{s-1}^{Q'_s} \bolds{1}_{\{V_{s-1}^{Q'_s} \leq h\}}
\leq h n.
\]
Therefore, if
\[
t_0 \stackrel{\triangle}{=}n \sqrt{5 h} \leq25 \sqrt{kmn},
\]
then
\[
h n \leq\frac{T^2}{5n}\qquad \mbox{ for any $T \geq t_0$},
\]
and so for such $T$ and any large enough $n$
\[
\sum_{s=1}^n V_{s-1}^{Q'_s} \bolds{1}_{\{V_{s-1}^{Q'_s} > h\}} \geq
\frac
{T^2}{4n} .
\]
By following the next arguments from the proof of Theorem~\ref
{thm-collisions-lower-bound},
it now follows that, as long as $T \geq t_0$,
\[
\mathbb{E}\operatorname{Col}_2(n) = \mathbb{E}\Biggl[\sum_{s=1}^t
X_{s-1}^{Q'_s}\Biggr] \geq\frac
{T^2}{8n}\bigl(1- O(n^{-4})\bigr) \geq\frac{T^2}{9n}.
\]
Similarly, using the argument as in the proof of Theorem~\ref
{thm-collisions-lower-bound}, which defines the stopping-time
$\tau_L$ and applies Proposition~\ref{prop-predictions} on the
sequence of increments given by
\[
\operatorname{Col}_2(t)-\operatorname{Col}_2(t-1) = N_{s-1}(J_s)
\bolds{1}_{A_s}(J_s),
\]
we deduce that, if $T \geq(t_0  \vee 30\sqrt{L n\log n})$ then
\[
\mathbb{P}\biggl(\operatorname{Col}_2(n) < \frac{T^2}{16n},\tau_L
> n\biggr) = O
(n^{-4}),
\]
as required.
\end{pf}

We next show how to infer the results regarding an unbounded maximal load
from Theorem~\ref{thm-collisions-gen}. For each integer $\ell
=0,1,2,\ldots,$ we define
the increasing sequence $(A_t)$ by
\[
A_t^{\ell} \stackrel{\triangle}{=}\{ i \in[n] \dvtx  N_t(i) \geq\ell\}.
\]
Further define
\[
T_\ell\stackrel{\triangle}{=}\sum_{s=1}^n Q_s (A_s^{\ell}),
\]
which is the expected number of balls that are placed in bins which
already hold at least $\ell$ balls.
The proof will follow from an inductive argument,
which bounds the value of $T_{\ell+1}$ in terms of $T_\ell$.

For some $L=L(n)$ to be specified later, our bounds will be meaningful
as long as the maximal load is at most $L$, and
%
\begin{equation}
\label{eq-T-ell-cond}
T_\ell\geq30\bigl( \sqrt{k m n}  \vee  \sqrt{L n\log n} \bigr).
\end{equation}
Using Theorem~\ref{thm-collisions-gen}, we will show that, if \eqref
{eq-T-ell-cond} holds then
%
\begin{equation}
\label{eq-T-ell-induction}
T_{\ell+1} \geq\frac{T_{\ell}^2}{20n L}.
\end{equation}
To this end, define
\[
R_\ell\stackrel{\triangle}{=}\sum_{i=1}^n \pmatrix{ N_n(i) - \ell\cr 2},
\]
that is, $R_\ell$ denotes the number of collisions between all pairs
of balls that were placed in a bin, that already held at least $\ell$ balls.

To infer \eqref{eq-T-ell-induction}, apply Theorem~\ref
{thm-collisions-gen} with respect to the subsets~$(A_s^\ell)$.
The assumption \eqref{eq-T-ell-cond} implies that, except with
probability $O(n^{-4})$, either the load is
at least~$L$, or
\[
R_\ell\geq\frac{T_\ell^2}{16n}.
\]
Notice that any ball that is placed in a bin, which contains at most
$L$ balls, can contribute at most $L$ collisions to the count of
$R_\ell$.
Therefore, if the maximal load is less than $L$, the following holds:
The number of balls placed in bins that already
contain at least $\ell$ balls, is at least
%
\begin{equation}
\label{eq-R-ell-bound}
R_\ell/ L \geq\frac{T_\ell^2}{16nL}\qquad \mbox{with probability
$1-O(n^{-4})$}.
\end{equation}
Recalling that $T_{\ell+1}$ is the expected number of such balls, we
infer that
\[
T_{\ell+1} \geq\bigl(1-O(n^{-4})\bigr) \frac{R_\ell}{L}
\geq\frac{T_\ell^2}{20nL},
\]
where the last inequality holds for large enough $n$ (with room to spare).


This establishes that \eqref{eq-T-ell-cond} implies \eqref
{eq-T-ell-induction}. Since by definition $T_0 = n$, we deduce
that the decreasing series $(T_0,T_1,\ldots)$ satisfies
\[
T_{\ell+1} \geq\frac{n}{(20 L)^{2^{\ell+1}-1}}\qquad  \mbox{if
$T_{\ell}$ satisfies \eqref{eq-T-ell-cond}}.
\]
Rearranging, it follows that, in particular, \eqref{eq-T-ell-cond} is
satisfied if
%
\begin{equation}\label{eq-suff-T-ell-cond}
30 \cdot(20L)^{2^\ell-1} \leq\sqrt{\frac{n}{km}}  \wedge  \sqrt
{\frac{n}{L\log n}}.
\end{equation}
It is now easy to verify that, for any fixed $\varepsilon> 0$, choosing
\[
L = \ell= (1-\varepsilon)\log_2\log\biggl(\frac{n}{km}\biggr)
\]
satisfies \eqref{eq-suff-T-ell-cond} for large enough $n$.
By \eqref{eq-R-ell-bound}, we can then infer that $R_\ell> 0$
with probability $1-O(n^{-4})$, hence the maximal load is at least
$\ell$.
This concludes the proof of Theorem~\ref{thm-matching-lower}. \qed

\section{Algorithms for perfect matching and constant load}\label
{sec:upper-bound}

In this section, we prove Theorem~\ref{thm-matching-upper} by
providing an algorithm that avoids collisions w.h.p. using only
$O(n/k)$ bits of memory,
which is the minimum possible by Theorem~\ref{thm-matching-lower}.
The case $k m = \Omega(n)$ of Theorem~\ref{thm-const-load} will then
follow from repeated applications of this algorithm.

\begin{tttt*}
\begin{enumerate}[ 1.]
\item For $\ell= \lfloor\frac{n}{\lfloor m/2\rfloor} \rfloor$,
partition the bins into contiguous blocks $B_1,\ldots,B_\ell$ each
comprising $\lfloor m/2\rfloor$ bins.
Ignore any remaining unused bins.
\item Set $d = \lceil\log_2(\frac{5}{C\delta}\log
n)\rceil$, and define the arrays $A_0,\ldots,A_{d-1}$:
\begin{itemize}
\item$A_j$ comprises $2^j$ contiguous blocks (a total of $\sim
2^{j-1}m$ bins).
\item For each contiguous (nonoverlapping) $4^j$-tuple of bins in
$A_j$, we keep a single bit that holds whether any of its bins is occupied.
\item All blocks currently or previously used are contiguous.
\end{itemize}
\item Repeat the following procedure until exhausting all rounds:
\begin{itemize}
\item Let $j$ be the minimal integer so that a bin of $A_j$, marked as
empty, appears in the current selection of $k$ bins. If no such $j$
exists, the algorithm announces failure.
\item Allocate the ball into this bin, and mark its $4^j$-tuple as occupied.
\item If the fraction of empty $4^j$-tuples remaining in $A_j$ just
dropped below $\delta/2$, relocate the array $A_j$ to a fresh
contiguous set of empty $2^j$ blocks (immediately beyond the last
allocated block).
If there are less than $2^j$ available new blocks, the algorithm fails.
\end{itemize}
\item Once $(1-\delta)n$ rounds are performed, the algorithm stops.
\end{enumerate}
\end{tttt*}

We proceed to verify the validity of the algorithm
in stages: First, we discuss a more basic version of the algorithm suited
for the case where $k m = \Omega(n \log n)$; then, we examine an
intermediate version
which extends the range of the parameters to $km\log m = \Omega(n\log
n)$; finally, we study the actual algorithm, which features the tight
requirement $km = \Omega(n)$.

Throughout the proof of the algorithm, assume that in each round
we are presented with $k$ independent uniform indices of bins, possibly
with repetitions.
Clearly, an upper bound for the maximal load in this relaxed model
translates into one for
the original model ($k$ choices without repetitions).

\subsection{Basic version of the algorithm}
We begin with a description and a proof of a simpler version of the
above algorithm, suited for the case where
%
\begin{equation}\label{eq-kn-nlogn}
k m \geq(3/\delta) n \log n.
\end{equation}
This version will serve as the base for the analysis.
For simplicity, assume first that $m \mid n$.

\begin{bbb*}
\begin{enumerate}
\item Let $B_1,\ldots,B_\ell$ be an arbitrary partition of the $n$
bins into $\ell\stackrel{\triangle}{=}n/m $ blocks,
each containing $m$ bins.
Put $r \stackrel{\triangle}{=}\lfloor(1-\delta) m \rfloor$.
\item Throughout stage $j\in[\ell]$, only the $m$ bins belonging to
$B_j$ are tracked.
At the beginning of the stage, all bins in the block are marked empty.
\item Stage $j$ comprises $r$ rounds, in each of which:
\begin{itemize}
\item The algorithm attempts to place a ball in an arbitrary empty bin
of $B_j$ if possible.
\item If no empty bin of $B_j$ is offered, the algorithm declares failure.
\end{itemize}
\item Once $(1-\delta)n$ rounds are performed, the algorithm stops.
\end{enumerate}
\end{bbb*}

To verify that this algorithm indeed produces a perfect allocation
w.h.p., examine a specific round of stage $j$, and condition on the event
that so far
the algorithm did not fail. In particular, its accounting of which bins
are occupied in $B_j$ is accurate,
and at least $m-r = (\delta-o(1))m$ bins in $B_j$ are still empty
[notice that by our
assumption $m =\Omega(\log n)$, and so $m\to\infty$ with $n$].

Let $\textsf{Miss}_j$ denote the event that
the next ball precludes all of the empty bins of $B_j$ in its $k$
choices, we have
%
\begin{equation}
\label{eq-alg-upper-prob}
\mathbb{P}(\textsf{Miss}_j) \leq\biggl(1-\frac{m-r}n\biggr)^k \leq
\mathrm
{e}^{-(\delta-o(1))({km}/n)} \leq n^{-3+o(1)},
\end{equation}
by assumption \eqref{eq-kn-nlogn}. A union bound over the $n$ rounds
now yields (with room to spare)
that the algorithm succeeds w.h.p.

The case where $m$ does not divide $n$ is treated similarly: Set $\ell
= \lfloor\frac{n}{\lfloor m/2 \rfloor} \rfloor$,
and\vspace*{1pt} partition the bins into blocks that now hold $\lfloor m/2\rfloor$
bins each,
except for the final block $B_\ell$ which would have between $\lfloor
m/2\rfloor$ and $m-1$ bins.
As before, in stage $j$ we attempt to allocate $\lfloor(1-\delta
)|B_j| \rfloor$ balls into $B_j$, while relying on the property
that $B_j$ has at least $(\delta-o(1))|B_j| \geq(\delta-o(1))m/2$
empty bins. This gives
\[
\mathbb{P}(\textsf{Miss}_j) \leq\mathrm{e}^{-(\delta-o(1))
({(km/2)}/n)} \leq
n^{-3/2+o(1)},
\]
as required.

\subsection{Intermediate version of the algorithm}
We now wish to adapt the above algorithm to the following case:
%
\begin{equation}\label{eq-km-n-m-poly}
k m \log_2 m \geq(20/\delta)\log(5/\delta) n \log n,\qquad  \log
^3 n \leq m \leq\frac{n}{\log n}.
\end{equation}
Notice that if $m \geq n^\varepsilon$, the above requirement is essentially
that
\[
k m = \Omega(n/\varepsilon).
\]
The full version of the algorithm will eliminate
this dependency on $\varepsilon$.

\begin{iii*}\vspace*{-9pt}
\begin{enumerate}
\item For $\ell= \lfloor\frac{n}{\lfloor m/2\rfloor} \rfloor$,
partition the bins into contiguous blocks $B_1,\ldots,B_\ell$ each\vspace*{1pt}
comprising $\lfloor m/2\rfloor$ bins.
Ignore any remaining unused bins.
\item Set $d = \lfloor\frac14 \log_2 m \rfloor$, and define the
arrays $A_0,\ldots,A_{d-1}$:
\begin{itemize}
\item$A_j$ is one of the blocks $B_1,\ldots,B_\ell$.
\item For each contiguous (non-overlapping) $2^j$-tuple of bins in
$A_j$, we keep a single bit that holds whether any of its bins is occupied.
\end{itemize}
\item Repeat the following procedure until exhausting all rounds:
\begin{itemize}
\item Let $j$ be the minimal integer so that a bin of $A_j$, marked as
empty, appears in the current selection of $k$ bins. If no such $j$
exists, the algorithm announces failure.
\item Allocate the ball into this bin, and mark its $2^j$-tuple as occupied.
\item If the fraction of empty $2^j$-tuples remaining in $A_j$ just
dropped below $\delta/2$, relocate the array $A_j$ to a fresh block
(immediately beyond the last allocated block).
If no such block is found, the algorithm fails.
\end{itemize}
\item Once $(1-\delta)n$ rounds are performed, the algorithm stops.
\end{enumerate}
\end{iii*}

Since the array $A_j$ contains $2^{-j}(m/2)$ different $2^j$-tuples,
the amount of memory required to maintain the status of all tuples is
\[
\frac{m}2 \sum_{j=0}^{d-1} 2^{-j} = (1- 2^{-d})m \leq m - m^{3/4}.
\]
In addition, we keep an index for each $A_j$, holding its position
among the $\ell$ blocks. By definition of $d$ and $\ell$, this
amounts to at most
\[
d \log_2\ell\leq(\log_2 n)^2 < m^{3/4}
\]
bits of memory, where the last inequality holds for any large $n$ by
\eqref{eq-km-n-m-poly}.

We first show that the algorithm does not fail to find a bin of $A_j$
marked as~empty. At any given point, each $A_j$ has a fraction of at
least $\delta/2$ bins marked~as empty. Hence, recalling \eqref
{eq-alg-upper-prob}, the probability of missing all the bins marked as~empty in $A_0,\ldots,A_{d-1}$ is at most
%
\begin{eqnarray}\label{eq-miss-Aj}
&&\exp\biggl[-\biggl(\frac\delta2-o(1)\biggr)\frac{km}{2n}
d\biggr]\nonumber\\
&&\qquad \leq
\exp\biggl[-\biggl(\frac\delta2-o(1)\biggr)\frac{10\log n}{\delta\log
_2 m} \log\biggl(\frac{20}\delta\biggr) \frac14\log_2 m
\biggr]\\
&&\qquad \leq n^{-\log(5/\delta)5/4-o(1)} < n^{-5/4},\nonumber
\end{eqnarray}
where the last inequality holds for large $n$.
Therefore, w.h.p. the algorithm never fails to find an array
$A_j$ with
an empty bin among the $k$ choices.

It remains to show that, whenever the algorithm relocates an array
$A_j$, there is always a fresh block available.

By the above analysis, the probability that a ball is allocated in
$A_j$ for $j\geq1$ at a given round is at most
\begin{eqnarray*}
\exp\biggl[-\biggl(\frac\delta2-o(1)\biggr)\frac{km/2}n j\biggr] &\leq&
\exp\biggl[-\biggl(\frac\delta2-o(1)\biggr)\frac{10\log n}{\delta\log
_2 m} \log\biggl(\frac{20}\delta\biggr) j \biggr]\\
&\leq&\exp\bigl(-3\log(5/\delta) j\bigr) \stackrel{\triangle}{=}p_j,
\end{eqnarray*}
where the last inequality holds for any sufficiently large $n$.

Let $N_j$ denote the number of balls that were allocated in blocks of
type $j$ throughout the run of the algorithm. Clearly, $N_j$ is
stochastically dominated by a binomial random variable $\operatorname
{Bin}(n,p_j)$.
Hence, known estimates for the binomial distribution (see, e.g., \cite
{AS}) imply that for all $j$,
\[
\mathbb{P}\bigl( N_j > n p_j + C \sqrt{n}\log n\bigr) \leq n^{-C}.
\]
The total number of blocks needed for $A_j$ is at most
\[
\bigg\lceil\frac{2^j N_j}{(1-\delta/2 )({m}/2)}\bigg\rceil,
\]
and hence the total number of blocks needed is w.h.p. at most
\begin{eqnarray*}
&&\Bigg\lceil\sum_{j=0}^{d-1}\frac{2^j (1-\delta)n p_j + C2^j\sqrt
{n}\log n}{(1-\delta/2 )({m}/2)} \Bigg\rceil\\
&&\qquad \leq
\sum_{j=0}^{d-1}\frac{2^j (1-\delta)n p_j}{(1-\delta/2 )({m}/2)} + O\biggl(\frac{n^{3/4}\log n}m\biggr).
\end{eqnarray*}
Since
\[
\sum_{j=1}^{d-1} 2^j p_j = \sum_{j=1}^{d-1}\exp\bigl( j \bigl(\log2 -
3\log(5/\delta)\bigr) \bigr) < 2\cdot2(\delta/5)^3 < \delta/5
\]
(with room to spare), the total number of blocks needed is w.h.p. at most
\[
\frac{(1+\delta/5)(1-\delta)n} {(1-\delta/2 )({m}/2)} +
O\biggl(\frac{n^{3/4}\log n}m\biggr) < \bigg\lfloor\frac{n}{\lfloor
m/2\rfloor} \bigg\rfloor= \ell
\]
for any sufficiently large $n$.

\subsection{Final version of the algorithm}

The main disadvantage in the intermediate version of the algorithm is that
the size of each $A_j$ was fixed at $m/2$ bins. Since the resolution
of each $A_j$ is in $2^j$-tuples, we are limited to at most $\log_2 m$
arrays. However, the probability of missing all the arrays $A_0,\ldots
,A_{d-1}$ has to compete with $n$, hence the requirement that $m$ would
be polynomial in $n$.

To remedy this, the algorithm uses arrays with increasing sizes, namely
$2^j$ blocks for $A_j$.
The resolution of each array is now in $4^j$-tuples, that is, tracking the
status of $A_j$ now requires at most $2^j \lfloor m/2\rfloor/ 4^j \leq
m/2^{j+1}$ bits.
Recalling that $d = \lceil\log_2(\frac{5}{C\delta}\log
n)\rceil$,
the number of memory bits required for all arrays is at most
%
\begin{equation}
\label{eq-mem-req}
\frac{m}2 \sum_{j=0}^{d-1} 2^{-j} = (1- 2^{-d})m \leq m - O(m/\log n).
\end{equation}
The following calculation shows that indeed there are sufficiently
many blocks to initially accommodate all the arrays:
\[
(2^d-1)\lfloor m/2\rfloor\leq\frac{5}{2C\delta}m \log n \leq\frac
{5km}{6C} = \frac56 n,
\]
where we used the assumptions $k \geq(3/\delta)\log n$ and $k m = C n$.

Each of the arrays comes along with a pointer to its starting block,
and the total number of memory bits required for this is at most
\[
d \log_2 (2n/m) \leq\bigl(\log_2 \log n + O(1)\bigr)\log_2 n = \bigl(1+o(1)\bigr)\log
_2 n \cdot\log_2 \log n.
\]
When $m = \Omega(\log^3 n)$, the space for these pointers clearly
fits among the $O(m/\log n)$ bits remaining according to \eqref{eq-mem-req}.
For smaller values of $m$, as before we can apply the algorithm for,
say, $m'=m/3$ (after tripling the constant $C_\delta$ to reflect this change),
thus earning $2m/3$ bits for the pointers (recall the requirement that
$m \geq\log_2 \log n \cdot\log_2 n$).

As final evidence that the choice of parameters for the algorithm is
valid, note\vspace*{1pt} that each $A_j$ indeed contains many $4^j$-tuples. It suffices
to check $A_{d-1}$, which indeed comprises about
\[
\bigl(1+o(1)\bigr)\frac{m}2 \frac{2^{d-1}}{4^{d-1}} = \bigl(1+o(1)\bigr)m/2^d =
\biggl(\frac{C\delta}5+o(1)\biggr)\frac{m}{\log n}=\Omega(\log\log n)
\]
$4^{d-1}$-tuples, where the last equality is by the assumption on the
order of $m$.

It remains to verify that the algorithm succeeds w.h.p. This will
follow from
the same argument as in the intermediate version of the algorithm. In
that version, each $A_j$
contained at least a fraction of $(\delta/2)$ empty bins, and $|A_j|$
was about $m/2$ for all $j$.
In the final version of the algorithm, each $A_j$ again contains at
least a fraction of $(\delta/2)$ empty bins,
but crucially, now $A_j$ contains $2^j$ bins. Thus, recalling \eqref
{eq-miss-Aj}, the probability to miss
$A_0,\ldots,A_{d-1}$ in a given round is now at most
\begin{eqnarray*}
\exp\Biggl[-\biggl(\frac\delta2-o(1)\biggr)\frac{km}{2n} \sum
_{j=0}^{d-1}2^j \Biggr]
&\leq&\exp\biggl(-\bigl(1-o(1)\bigr)\frac{C\delta}4 (2^d-1) \biggr)\\
& =& n^{-5/4-o(1)},
\end{eqnarray*}
where the last inequality is by the definition of $d$.
A union bound over the $n$ rounds gives that, w.h.p., an array $A_j$
with an empty bin is found for every ball.

To see that w.h.p. there are always sufficiently many available fresh
blocks to relocate an array,
one essentially repeats the argument from the intermediate version of
the algorithm. That is,
we again examine the probability that a ball is allocated in $A_j$, to
obtain that this time
\[
p_j = \exp\biggl(-\bigl(1-o(1)\bigr)\frac{C\delta}4 (2^j-1)\biggr).
\]
A choice of $C \geq K (1/\delta)\log(1/\delta)$ with
some suitably large $K > 0$ would give
\[
\sum_{j=1}^{d-1}4^j p_j < \delta/ 5,
\]
and the rest of that argument unchanged now implies that the algorithm
never runs out of fresh blocks w.h.p.

This completes the proof of Theorem~\ref{thm-matching-upper}.

\subsection[Proof of upper bound in Theorem 1]{Proof of upper bound in Theorem~\protect\ref{thm-const-load}}
We now wish to apply the algorithm from Theorem~\ref
{thm-matching-upper} in order to obtain a
constant load in the case where $k m \geq c n$ for some $c> 0$. To
achieve this, consider the perfect
matching algorithm for, say, $\delta=\frac12$, and let $C_\delta$ be
the constant that appears in Theorem~\ref{thm-matching-upper}.
Next, join every consecutive $\lceil C_\delta/c\rceil$-tuple of bins
together and write $n'$ for the number of such tuples.
As $km \geq C n'$, we may apply the perfect-matching algorithm for
$n'/2$ balls with respect to the $n'$ tuples of bins, keeping
in mind that the algorithm is valid also for the model of repetitions.
This gives a perfect matching w.h.p., and repeating this process
gives a
total load of at most $2C_\delta/c = O(1)$ for all $n$ balls. \qed

\eject
\section{Improved lower bounds for poly-logarithmic choices}\label
{sec:large-q}
%
%
\subsection[Proof of Theorem 2]{Proof of Theorem~\protect\ref{thm-large-load}}
Our proof of this case is an extension of the proof of Theorem~\ref
{thm-matching-lower}. We now wish to estimate the number of
$q$-collisions for general $q$:
\[
\operatorname{Col}_q(t) \stackrel{\triangle}{=}\sum_{i=1}^n \pmatrix{ N_t(i)\cr q}.
\]
The analysis hinges on a recursion on $q$, for which we need to achieve
bounds on a generalized quantity, a linear function of the
$q$-collisions vector:
%
\begin{eqnarray}
\qquad X_t^{f;q} &\stackrel{\triangle}{=}&\sum_{s_1<\cdots< s_{q} \leq t}
\sum_i f(i) \bolds{1}_{\{
J_{s_1}=i\}}\cdots\bolds{1}_{\{J_{s_{q}}=i\}} = \sum_i f(i) \pmatrix{ N_{t}(i)\cr q},\label{eq-Xt-def}\\
V_t^{f;q} &\stackrel{\triangle}{=}&\sum_{s_1<\cdots< s_{q} \leq t}
\sum_i f(i)
Q_{s_1}(i)\cdots Q_{s_{q}}(i) \label{eq-Vt-def}.
\end{eqnarray}
Our objective is to obtain lower bounds for $X_t^{f;q}$ with $f \equiv
1$, as clearly \mbox{$\operatorname{Col}_q(t) = X_t^{\bolds{1};q}$}.
Notice that the parameters $X_t^\nu,V_t^\nu$ from Section~\ref
{sec:lower-bound} are exactly\vspace*{1pt} $X_t^{\nu;1},V_t^{\nu;1}$ defined above.
There, $\nu$ was a strategy, whereas now
our $f$ will be the product of different strategies. This fact will
allow us to formulate
a recursion relation between the $V_t^{f;q}$'s and an approximate
recursion for the $X_t^{f;q}$.
We achieve this using the next lemma, where here in and throughout the
proof we let
%
\begin{equation}
\label{eq-L-def}
L \stackrel{\triangle}{=}\log(n/m)
\end{equation}
denote a maximal load we do not expect to reach (except if the
algorithm is far from optimal).
We further define
\[
\Gamma\stackrel{\triangle}{=}\Biggl\{ \prod_{i=1}^L f_i
 \dvtx   f_i\in\{
\bolds{1},\mu_1,\ldots,\mu_{2^m}\}\mbox{ for all $i$}\Biggr\}
\]
to be the set of all point-wise products of at most $L$ strategies from
the pool.

\begin{lemma}\label{lem-q-collisions-induction}
Either the maximal load exceeds $L$, or
the following holds for all $q < L$, every $t \leq n/k$ and every $f
\in\Gamma$,
except with probability $\mathrm{e}^{-3mL}$:
%
\begin{equation}
\label{eq-Vt-large-Xt-large}
\mbox{If }
V_t^{f;q} \geq100 \frac{(3L)^{q+1}}{q!} m \|f\|_\infty\qquad  \mbox{then }
X_t^{f;q} \geq3^{-q} V_t^{f;q}.
\end{equation}
\end{lemma}

\begin{pf}
The key property of the quantities $V_t^{f;q}$, which justified the
inclusion of the inner products with $f$, is the following
recursion relation, whose validity readily follows from definition
\eqref{eq-Vt-def}:
%
\begin{equation}
\label{eq-Vt-recursion}V_t^{f;q+1} = \sum_{s < t}
V_{s}^{(Q_{s+1}\cdot f) ; q}\qquad  \mbox{ for any $q\geq1$ and any $t$}.
\end{equation}
%
We now wish to write a similar recursion for the variables $X_t^{f;q}$.
As opposed to the variables $V_t^{f;q}$,
which satisfied the above recursion combinatorially, here the recursion
will only be stochastic. Notice that
\begin{eqnarray*}
X_{t+1}^{f;q+1} - X_t^{f;q+1}
&=& f(J_{t+1})\biggl(\pmatrix{ N_{t}(J_{t+1})+1\cr q+1}-\pmatrix{ N_{t}(J_{t+1})\cr q+1}\biggr)\\
&= &f(J_{t+1})\pmatrix{ N_{t}(J_{t+1})\cr q},
\end{eqnarray*}
and hence
\[
\mathbb{E}[X_{t+1}^{f;q+1}-X_t^{f;q+1} \mid\mathcal{F}_t
] = \sum_i
Q_{t+1}(i)f(i)\pmatrix{ N_{t}(i)\cr q} = X_{t}^{(Q_{t+1}\cdot f);q}.
\]
%
We may therefore apply Proposition~\ref{prop-predictions} as follows:
\begin{itemize}
\item The sequence of increments we consider is
$(X_{t+1}^{f;q+1}-X_t^{f;q+1})$ (that results in a telescopic sum).
\item The sequence of conditional expectations is $(X_t^{(Q_{t+1}\cdot
f); q})$.
\item The bound on the increment is $M = \|f\|_\infty{L\choose q}$,
where $L$ is an upper bound for the maximal load (if we encounter
a load of $L$, we stop the process).
\end{itemize}
This implies that
\begin{eqnarray*}
&&\mathbb{P}\biggl(\exists t \dvtx  \biggl\{X_t^{f;q+1} \leq\frac12 \sum
_{s<t}X_s^{(Q_{s+1}\cdot f);q},\sum_{s<t}X_s^{(Q_{s+1}\cdot
f);q} \geq100 m L \| f\|_\infty\pmatrix{ L\cr q}\biggr\}\biggr) \\
&&\qquad \leq\exp\biggl(-\frac{100 m L \| f\|_\infty{ L\choose q}}{20\|f\|
_\infty{ L\choose q}}+2\biggr) = O(\exp(-5mL)).
\end{eqnarray*}
As a result, the above event does not occur for any $f \in\Gamma$
(since there are at most $2^{mL}$ such functions)
except with probability $\mathrm{e}^{-4mL}$. Therefore, setting
\[
h_{f;q} \stackrel{\triangle}{=}100 \frac{(3L)^{q+1}}{q!}m\|f\|
_\infty,
\]
we have that, except with probability $\mathrm{e}^{-4mL}$,
%
\begin{equation}
\label{eq-Xt-large-sum-Xt-large}
\mbox{if }
\sum_{s<t}X_s^{(Q_{s+1}\cdot f);q} > 3^{-q}h_{f;q}\qquad   \mbox{then
}
X_t^{f;q+1} \geq\frac12 \sum_{s<t}X_s^{(Q_{s+1}\cdot f);q}.
\end{equation}

We now proceed to prove \eqref{eq-Vt-large-Xt-large} by induction on
$q$. For $q=1$, notice that
\begin{eqnarray*}
X_t^{f;1} &=& \sum_{s \leq t} \sum_i f(i) \bolds{1}_{\{J_{s}=i\}}=
\sum_i
f(i) N_{t}(i),\\
V_t^{f;1} &=& \sum_{s \leq t} \sum_i f(i) Q_{s}(i).
\end{eqnarray*}
Furthermore, as the definition of $X_t^{f;q}$ also applies to the case
$q=0$, we obtain that
\begin{eqnarray*}
X_t^{f;0} &=& \sum_i f(i)\qquad \mbox{and so}\\
V_t^{f;1} &=& \sum_{s < t} \sum_i Q_{s+1}(i)f(i) = \sum
_{s<t}X_s^{(Q_{s+1}\cdot f);0}.
\end{eqnarray*}
Hence, combining the assumption $V_t^{f;1} \geq100 (3L)^2 m \|f\|
_\infty= h_{f;1} $
with statement \eqref{eq-Xt-large-sum-Xt-large} yields that $X_t^{f;1}
\geq\frac12 V_t^{f;1} \geq\frac13 V_t^{f;1}$, except with
probability $\mathrm{e}^{-4mL}$.

It remains to establish the induction step. The induction hypothesis
for $q$ states that whenever $V_t^{f;q} \geq h_{f;q}$ we also have
$X_t^{f;q}\geq3^{-q} V_t^{f;q}$
except with probability~$\mathrm{e}^{-4mL}$. Therefore,
%
\begin{eqnarray}\label{eq-sum-Xs-1}
\sum_{s<t}X_s^{(Q_{s+1}\cdot f) ; q} &\geq&
3^{-q} \sum_{s<t} V_s^{(Q_{s+1}\cdot f) ; q} \cdot\bolds{1}_{\{
V_s^{(Q_{s+1}\cdot f) ; q} > h_{(Q_{s+1}\cdot f); q}\}} \nonumber\\
&\geq&3^{-q} \biggl(\sum_{s<t} V_s^{(Q_{s+1}\cdot f) ; q} - t \cdot
h_{(Q_{s+1}\cdot f); q}\biggr)\\
&\geq&3^{-q} \biggl(V_t^{f;q+1} - t \cdot100 \frac{(3L)^{q+1}}{q!}m \|
Q_{s+1}\cdot f\|_\infty\biggr),\nonumber
\end{eqnarray}
where in the last inequality we applied the recursion relation \eqref
{eq-Vt-recursion}.
Recalling that $Q_{s+1}$ is a strategy, the following holds for all $t
\leq n/k$:
\[
t\|Q_{s+1}\cdot f\|_\infty\leq t\|Q_{s+1}\|_\infty\|f\|_\infty\leq
t\frac{k}n\|f\|_\infty\leq\|f\|_\infty.
\]
Plugging this into \eqref{eq-sum-Xs-1}, we obtain that for all $t \leq n/k$,
\begin{eqnarray}\label{eq-sum-Xs-2}
\sum_{s<t}X_s^{(Q_{s+1}\cdot f) ; q} &\geq&3^{-q} \biggl(V_t^{f;q+1} -
100 \frac{(3L)^{q+1}}{q!}m \|f\|_\infty\biggr)\nonumber\\[-8pt]\\[-8pt]
&=& 3^{-q} (V_t^{f;q+1} - h_{f;q}).\nonumber
\end{eqnarray}
Now, if $V_t^{f;q+1} \geq h_{f;q+1} = 100 \frac{(3L)^{q+2}}{(q+1)!}m\|
f\|_\infty$, then in particular
\[
V_t^{f;q+1} \geq3\cdot100 \frac{(3L)^{q+1}}{q!}m\|f\|_\infty=
3h_{f;q} \qquad  \mbox{for all $q \leq L-1$.}
\]
Thus, under this assumption, \eqref{eq-sum-Xs-2} takes the following form:
\[
\sum_{s<t}X_s^{(Q_{s+1}\cdot f) ; q} \geq3^{-q} \cdot\frac23
V_t^{f;q+1} \geq3^{-q}\cdot2 h_{f;q} .
\]
Since this satisfies the condition of \eqref{eq-Xt-large-sum-Xt-large}
(where we actually only needed a lower
bound of $3^{-q} h_{f;q}$), we obtain that except with probability
$\mathrm{e}^{-4mL}$,
\[
X_t^{f;q+1} \geq\frac12 \sum_{s<t}X_s^{(Q_{s+1}\cdot f) ; q} \geq
3^{-(q+1)} V_t^{f;q+1},
\]
completing the induction step.

Summing the error probabilities over the induction steps for every $q <
L$ concludes the proof of the lemma.
\end{pf}

It remains to apply the above lemma to deduce the maximal load of
$\Omega(\frac{\log n}{\log\log n})$ for $k=\operatorname{polylog}(n)$.
Recalling that $m \leq n^{1-\delta}$ for some fixed $\delta> 0$, let
$0<\varepsilon< \delta/2$ and choose the following parameters:
\[
q = (1-\varepsilon)\frac{\log(n/m)}{\log k + \log\log(n/m)},\qquad f =
\bolds{1},t = n/k.
\]
Lemma~\ref{lem-q-collisions-induction} now gives that, either the
maximal load exceeds $L = \log(n/m)$, or w.h.p. the following statements holds:
%
\begin{equation}
\label{eq-Vn/k-Xn/k-implication}
\mbox{If }
V_{n/k}^{\bolds{1};q} \geq100 \frac{(3L)^{q+1}}{q!} m \qquad  \mbox{then
}
X_{n/k}^{\bolds{1};q} \geq3^{-q} V_{n/k}^{\bolds{1};q}.
\end{equation}
Notice that for the above value of $q$, we have $3^{-q} =
(n/m)^{o(1)}$, and therefore, showing that the condition of \eqref
{eq-Vn/k-Xn/k-implication}
is satisfied and that $V_{n/k}^{\bolds{1};q} \geq
(n/m)^{\varepsilon/2}$ would immediately imply that the maximal load
exceeds $q$ w.h.p.

The following lemma, which
provides a lower bound on $V_t^{\bolds{1};q}$, is thus the final ingredient
required for the proof of the theorem:

\begin{lemma}\label{lem-Vt-packing}
For all $t$, $k$ and $q$, all $Q_1,\ldots,Q_t$ and any fixed $\alpha>
0$ we have
\[
V_{t}^{\bolds{1};q} \geq\frac{(t-(1+\alpha)kq
)^q}{\mathrm
{e}^{q/(2\alpha)}n^{q-1}q!}.
\]
\end{lemma}

\begin{pf}
Recall that
\begin{eqnarray*}
V_t^{\bolds{1};q} &=& \sum_{s_1<\cdots<s_q\leq t} \sum_{i=1}^n
(Q_{s_1}\cdots Q_{s_q})(i)\nonumber\\
&=&\frac1{q!}\sum_{i=1}^n \sum_{s_1\leq t} Q_{s_1}(i) \mathop{\sum
_{s_2\leq t}}_{s_2\neq s_1} Q_{s_2}(i) \cdots\mathop{\sum_{s_q \leq
t}}_{s_q\notin\{s_1,\ldots,s_{q-1}\}} Q_{s_q}(i).
\end{eqnarray*}
Defining
\[
r_i \stackrel{\triangle}{=}\sum_{s \leq t} Q_s(i),
\]
and recalling that $\|Q_{s}\|_\infty\leq k/n$ for all $s$, it follows
that for all $i$ and $j\geq1$,
\[
\mathop{\sum_{s_j\leq t}}_{s_j\neq\{s_1,\ldots,s_{j-1}\}}
Q_{s_j}(i) \geq r_i - (j-1)k/n.
\]
Consequently,
\begin{eqnarray}\label{eq-Vt1q-bound-1}
V_t^{\bolds{1};q} &\geq&\frac1{q!}\sum_{i=1}^n \prod
_{j=1}^{q} \bigl(r_i - (j-1) k/n\bigr)\nonumber\\[-8pt]\\[-8pt]
&\geq&\frac1{q!}\sum_{i=1}^n \prod_{j=1}^{q} \bigl(r_i - (j-1)
k/n\bigr) \bolds{1}_{\{r_i > (1+\alpha){(k(q-1)}/n)\}}.\nonumber
\end{eqnarray}
Next, notice that for all $1 \leq j\leq q-1$,
\begin{eqnarray*}
1-\frac{j}{(1+\alpha)(q-1)} &\geq&\exp\biggl[ -\frac{j}{(1+\alpha
)(q-1)} \Big/ \biggl(1- \frac{j}{(1+\alpha)(q-1)}\biggr) \biggr]\\
&\geq&\exp\biggl[ -\frac{j}{\alpha(q-1)}\biggr].
\end{eqnarray*}
Thus, in case $r_i > (1+\alpha)\frac{k(q-1)}n$ we have the following
for all $1 \leq j \leq q$:
\[
r_i - \frac{(j-1) k}n > r_i\biggl( 1 - \frac{j-1}{(1+\alpha
)(q-1)}\biggr) \geq r_i \exp\biggl[ -\frac{j-1}{\alpha(q-1)}\biggr].
\]
Combining this with \eqref{eq-Vt1q-bound-1}, we deduce that
\begin{eqnarray*}
V_t^{\bolds{1};q} &\geq&\frac1{q!}\sum_{i=1}^n \prod
_{j=1}^{q} r_i \exp\biggl[ -\frac{j-1}{\alpha(q-1)}\biggr]
\bolds{1}_{\{r_i > (1+\alpha)(k(q-1)/n)\}}\\
&=& \frac1{q!}\sum_{i=1}^n \bigl( \mathrm{e}^{-1/(2\alpha)} r_i
\bolds{1}_{\{r_i > (1+\alpha)({k(q-1)}/n)\}}\bigr)^q.
\end{eqnarray*}
Applying Cauchy--Schwarz, we infer that
\begin{eqnarray*}
V_t^{\bolds{1};q} &\geq&\frac{n}{q!} \biggl(\frac{\mathrm
{e}^{-1/(2\alpha)} \sum_i r_i
\bolds{1}_{\{r_i > (1+\alpha){{k(q-1)}/n}\}}}n \biggr)^q \\
&=&
\frac{(\sum_i r_i \bolds{1}_{\{r_i > (1+\alpha)
({k(q-1)}/n)\}}
)^q}{
\mathrm{e}^{q/(2\alpha)}n^{q-1}q!}.
\end{eqnarray*}
The proof of the lemma now follows from noticing that
\[
\sum_i r_i \bolds{1}_{\{r_i \leq(1+\alpha)({k(q-1)}/n)\}} \leq
(1+\alpha) k(q-1) < (1+\alpha)k q,
\]
whereas $\sum_i r_i = \sum_{s \leq t} \sum_i Q_s(i) = t$.
\end{pf}

To complete the proof using Lemma~\ref{lem-Vt-packing}, apply this
lemma for $\alpha=1$, $t = n/k$, and $kq=n^{o(1)}$, giving that
\[
V_{n/k}^{1;q} \geq\frac{n}{((\mathrm{e}^{1/2}-o(1)
)k)^{q} q!} \geq(2k)^{-q} n / q!,
\]
where the last inequality holds for any sufficiently large $n$. Consequently,
\[
\frac{100(3L)^{q+1}m/q!}{V_{n/k}^{\bolds{1};q}} \leq100(3L)^{q+1} (2k)^q
\frac{m}n \leq
100(6kL)^{q+1} \frac{m}n .
\]
Since our choice of $q$ is such that
\[
(6kL)^{q+1} = \mathrm{e}^{(1+o(1))q (\log L + \log k)} =
(n/m)^{1-\varepsilon-o(1)},
\]
we have that
\[
\frac{100(3L)^{q+1}m/q!}{V_{n/k}^{\bolds{1};q}} \leq
(n/m)^{-\varepsilon+o(1)}.
\]
This implies both that $V_{n/k}^{\bolds{1};q} \geq(n/m)^{\varepsilon
/2}$ for
any large $n$ (recall that $L > q$), and that the condition of
\eqref{eq-Vn/k-Xn/k-implication} is satisfied for any large $n$.
Altogether, the
maximal load is w.h.p. at least $q$, concluding the proof of
Theorem~\ref{thm-large-load}. \qed

\subsection{A corollary for nonadaptive algorithms}
We end this section with a corollary of Theorem~\ref{thm-large-load}
for the case
of non-adaptive algorithms, that is, the strategies $Q_1,\ldots,Q_n$ are fixed
ahead of time. Namely, we show that for $k = O(n\frac{\log\log
n}{\log n})$ the optimal maximal
load is w.h.p. $\Theta(\frac{\log n}{\log\log n})$, that is, of the
same order as the one for $k=1$.
Theorem~\ref{thm-nonadaptive}, whose proof appears in Section~\ref
{sec:nonadaptive}, includes
a different approach that proves this result more directly.

\begin{corollary}
Consider the allocation problem of $n$ balls
into $n$ bins, where each ball has $k$ independent uniform choices.
If $k \leq C n\frac{\log\log n}{\log n}$, then any nonadaptive algorithm
w.h.p. has a maximal-load of at least $\frac{1-o(1)}{C \vee
1}\cdot
\frac{\log n}{\log\log n}$.
In particular, if $k \leq n\frac{\log\log n}{\log n}$ then the load is
at least $(1-o(1))\frac{\log n}{\log\log n}$ w.h.p.
\end{corollary}

\begin{pf}
Let $Q_1,\ldots,Q_n$ be the optimal sequence of strategies for the
problem. Using definitions \eqref{eq-Xt-def} and \eqref{eq-Vt-def}
with $f \equiv1$, we have the following for all $q$:
\begin{eqnarray*}
X_t^{\bolds{1};q} &=& \sum_i \pmatrix{ N_t(i)\cr q} = \operatorname
{Col}_q(t)\quad \mbox{and
}\\
V_t^{\bolds{1};q} &=& \mathbb{E}X_t^{\bolds{1};q}.
\end{eqnarray*}
Fix $0 < \varepsilon< \frac12$. Applying Lemma~\ref{lem-Vt-packing}
with $t=n$ and $\alpha= \varepsilon/[2(1-\varepsilon)]$,
%
\begin{equation}
\label{eq-V1n-lower-bound-eps}
V_n^{\bolds{1};q} \geq\frac{(n-(1+\alpha)kq)^q}{\mathrm
{e}^{q/(2\alpha
)} n^{q-1}q!} =
\biggl(1-\frac{(2-\varepsilon)kq}{2(1-\varepsilon)n}\biggr)^q \cdot
\frac
{n}{\mathrm{e}^{(1-\varepsilon)q/\varepsilon}q!}.
\end{equation}
Recalling that $k \leq C n \frac{\log\log n}{\log n}$ for some fixed
$C > 0$, set
%
\begin{equation}
\label{eq-q-value-non-adaptive}
q = \frac{1-\varepsilon}{C \vee  1}\cdot\frac{\log n}{\log\log n}.
\end{equation}
This choice has $kq/n \leq1-\varepsilon$ and $q! \leq
n^{1-\varepsilon
+o(1)}$. Combined with \eqref{eq-V1n-lower-bound-eps},
%
\begin{eqnarray}\label{eq-non-adapt-E-q-col}
\mathbb{E}\operatorname{Col}_q(n) &=& V_n^{\bolds{1};q}\nonumber\\
 & \geq&\exp
\biggl(-\frac{(2-\varepsilon
)kq^2}{2(1-\varepsilon)n} \Big /  \biggl(1-\frac{(2-\varepsilon
)kq}{2(1-\varepsilon)n}\biggr)\biggr)
\frac{n}{\mathrm{e}^{(1-\varepsilon)q/\varepsilon} q!} \\
& \geq&\exp\biggl(-\frac{2-\varepsilon}{\varepsilon}q\biggr) \frac
{n}{n^{1-\varepsilon+ o(1)}} = n^{\varepsilon- o(1)} .\nonumber
\end{eqnarray}
To translate the number of $q$-collisions to the number of bins with
load $q$, consider
the case where for some bin $j$ we have $\sum_{s=1}^n Q_s(j) \geq100
\log n$. Proposition~\ref{prop-predictions} (applied
to the Bernoulli variables $\bolds{1}_{\{J_s = j\}}$) then implies
that $N_n(j) \geq50\log n$ except with probability $O(n^{-5})$, and in
particular the maximal load exceeds $q$ w.h.p.
We may therefore assume from this point on that $\sum_{s=1}^n Q_s(j)
\leq100 \log n$ for all~$j$.

Set $L \stackrel{\triangle}{=}150 \log n$. Clearly, upon increasing
$Q_s(j)$ for some $1
\leq s\leq n$, the load in bin $j$ will stochastically
dominate the original one. Thus, for any integer $r \geq1$ we may
increase $\sum_{s=1}^n Q_s(j)$ to $\frac23 r L$, and by
Proposition~\ref{prop-predictions} obtain that
\mbox{$N_n(j) \leq r L$} except with probability $O(\exp(-rL/30))$. Defining
\[
A_r \stackrel{\triangle}{=}\Bigl\{ rL \leq\max_{1\leq j \leq n}
N_n(j) < (r+1)L\Bigr\}\qquad
\mbox{(for $r=0,1,\ldots$)},
\]
we in particular get $\mathbb{P}(A_r) = O(n\exp(-rL/30))$. However, clearly
on this event $\operatorname{Col}_q(n) \leq n {(r+1)L\choose q}$,
and since $n^2{(r+1)L\choose q} \leq O(\exp(rL/50))$, we have
\[
\mathbb{E}[\operatorname{Col}_q(n) \mid\overline{A_0}
]\mathbb{P}(\overline
{A_0}) \leq\sum_{r\geq1} O\bigl(\exp(-rL/100)\bigr) = O(n^{-3/2}) = o(1).
\]
Thus, by \eqref{eq-non-adapt-E-q-col}, we have $\mathbb
{E}[\operatorname{Col}_q(n) \mid
A_0] \geq n^{\varepsilon-o(1)}$.
Finally, since any given bin can contribute at most ${L\choose q} =
n^{o(1)}$ collisions to $\operatorname{Col}_q(n)$ given $A_0$,
\[
\mathbb{E}\Biggl[\sum_{j=1}^n \bolds{1}_{\{N_n(j) \geq q\}}\Biggr]
\geq\mathbb{E}
\biggl[\operatorname{Col}_q(n)\Big/\pmatrix{ L\cr q} \Big| A_0\biggr]
\mathbb{P}(A_0) =
n^{\varepsilon-o(1)}.
\]

As demonstrated in the next section (see Lemma~\ref
{lem-negative-correlation}), one can now use the fact that the events
$\{N_n(j) \geq q\}$ are negatively
correlated to establish concentration for the variable $\sum_{j=1}^n
\bolds{1}_{\{N_n(j) \geq q\}}$.
Altogether, we deduce that the maximal load w.h.p. exceeds $q$,
as required.
\end{pf}

\section{Tight bounds for nonadaptive allocations}\label{sec:nonadaptive}

In this section, we present the proof of Theorem~\ref{thm-nonadaptive}.
Throughout the proof we assume, whenever this is needed, that
$n$ is sufficiently large. To simplify the presentation, we omit all
floor and ceiling signs whenever these are not crucial.
We need the following lemma.

\begin{lemma}
\label{l61}
Let $p_1,p_2, \ldots,p_n$ be reals satisfying $0 \leq p_i \leq
\frac{\log\log n}{\log n}$ for all~$i$, such that $\sum_{i=1}^n p_i
\geq
1-\varepsilon$, where $\varepsilon=\varepsilon(n) \in[0,1]$. Let
$X_1, X_2,
\ldots,X_n$
be independent indicator random variables, where
$\mathbb{P}(X_i=1)=p_i$ for all $i$, and put $X=\sum_{i=1}^n X_i$. Then
\[
\mathbb{P}\biggl(X \geq(1-\varepsilon) \frac{\log n}{\log\log n}
\biggr) \geq
\frac{1}{n^{1-\varepsilon}}.
\]
\end{lemma}

\begin{pf}
Without loss of generality, assume that $p_1 \geq p_2 \geq\cdots
\geq p_n$. Define a family
of $k$ pairwise disjoint blocks $B_1,B_2, \ldots,B_k \subset\{1, 2,
\ldots,n\}$, where
$k \geq(1-\varepsilon) \frac{\log n}{\log\log n} $
so that for each $i$, $1 \leq i \leq k$,
\[
\frac{2}{\log n} \leq\sum_{j \in B_i} p_j \leq\frac{\log\log
n}{\log n}.
\]
This can be easily done greedily; the first block consists of
the indices $1,2, \ldots,r$ where $r$ is the smallest integer
so that $\sum_{j=1}^r p_j \geq\frac{2}{\log n}$. Note that it is
possible that $r=1$, and that since the sequence $p_j$ is monotone
decreasing, $\sum_{j=1}^r p_j \leq\frac{\log\log n}{\log n}$.
Assuming we have already partitioned the indices $\{1,\ldots
,r\}$ into blocks, and assuming we still do not have $(1-\varepsilon)
\frac{\log n}{\log\log n}$ blocks, let the next block
be $\{r+1, \ldots,s\}$ with $s$ being the smallest integer
exceeding $r$ so that \mbox{$\sum_{j=r+1}^s p_j \geq
\frac{2}{\log n}$}. Note that if $p_{r+1} \geq\frac{2}{\log n}$
then $s=r+1$, that is, the block consists of a single element,
and otherwise $\sum_{j=r+1}^s p_j < \frac{4}{\log n} <
\frac{\log\log n}{\log n}$. Thus,
in any case the sum above is at least $\frac{2}{\log n}$ and at
most $\frac{\log\log n}{\log n}$. Since the total sum of the reals
$p_j$ is at least $1-\varepsilon$ this process does not terminate
before generating $k \geq(1-\varepsilon) \frac{\log n}{\log\log n}$
blocks, as needed.

Fix a family of $k=(1-\varepsilon) \frac{\log n}{\log\log n}$ blocks
as above. Note that for each fixed block $B_i$ in the family, the
probability that $\sum_{j \in B_i} X_j \geq1$ is at least
\[
\sum_{j \in B_i} p_j -\sum_{j,q \in B_i, j<q} p_j p_q
\geq\sum_{j \in B_i}p_j -\frac{1}{2} \biggl(\sum_{j \in B_i} p_j\biggr)^2
\geq\frac{2}{\log n} -\frac{2}{\log^2 n} > \frac{1}{\log n}.
\]
It thus follows that the probability that for each of the $k$
blocks $B_i$ in the family $\sum_{j \in B_i} X_j \geq1$ is
at least $(\frac{1}{\log n})^k=\frac{1}{n^{1-\varepsilon}}$,
completing the proof of the lemma.
\end{pf}

\begin{pf*}{Proof of Theorem~\ref{thm-nonadaptive} [\normalfont{Part~(i)}]}
As before, our framework is the relaxed model where there are strategies
$Q_1,Q_2,\ldots,Q_n$, where $Q_t$ is the distribution of the bin
to be selected for ball number $t$, satisfying $\|Q_t\|_\infty\leq k/n$.
However, since now we consider nonadaptive algorithms, the strategies
are no longer random variables, but rather a predetermined sequence.
We therefore let $P = (p_{i t})$ denote the $n\times n$ matrix of probabilities,
where $p_{i t}$ is the probability that the ball at time $t$ would be
placed in
bin $i$.
Clearly,
\[
0 \leq p_{it} \leq k/n =\frac{\log\log n}{\log n}\qquad \mbox{for all
$i$ and $t$}
\]
and
\[
\sum_{1 \leq i \leq n} p_{it}=1\qquad \mbox{for all $t$}.
\]
The sum of entries of each column
of the $n$ by $n$ matrix
$p_{it}$ is $1$, and hence the total sum of its entries is $n$.
If it contains a row $i$ so that the sum of entries in this row
is at least, say, $\log n$, then the expected number of balls in
bin number $i$ by the end of the process is $\sum_{t=1}^n p_{it}
\geq\log n$. As the variance is
\[
\sum_{t+1}^n p_{it} (1-p_{it})
\leq\sum_{t=1}^n p_{it},
\]
it follows by Chebyshev's inequality (or by Hoeffding's inequality)
that with high probability the actual number of balls placed in bin
number $i$ exceeds $\frac{\log n}{2} > \frac{\log n}{\log\log
n}$, showing that in this case the desired result holds.

We thus assume that the sum of entries in each row is at most $\log
n$. As the average sum in a row is $1$, there is a row whose total
sum is at least $1$. Omit this row, and note that since its total
sum is at most $\log n$, the sum of all remaining entries of the
matrix is still at least $n-\log n$, and hence the average sum
of a row in it is at least $\frac{n-\log n}{n-1} > 1-\frac{\log
n}{n}$. Therefore, there is another row of total sum at least this
quantity. Omitting this row and proceeding in this manner we can
define a set of rows so that the sum in each of them is large. Note
that as long as we defined at most $\frac{n}{\log^2 n}$ rows, the
total sum of the remaining elements of the matrix is still at least
$n-\frac{n}{\log n}$, and hence there is another row of total sum
at least $1-\frac{1}{\log n}$. We have thus shown that there is a
set $I$ of $\frac{n}{\log^2 n}$ rows such that
\[
\sum_{t=1}^n p_{it} \geq1-\frac{1}{\log n}\qquad \mbox{for each $i \in I$}.
\]
For each $i\in I$, let $A_i$ denote the event
\[
A_i \stackrel{\triangle}{=}\biggl( \mbox{There are at most }
\biggl(\frac{\log
n}{\log\log n}-4\biggr)\mbox{ balls in bin $i$ }\biggr).
\]
Applying Lemma \ref{l61} with $\varepsilon=\frac{4 \log\log n}{\log
n}$, we get
\[
\mathbb{P}(A_i) \leq1- \frac{ \log^4 n}{ n}\qquad \mbox{for each $i\in I$}.
\]
We will next show that, as the events $A_i$ are negatively correlated,
the probability that all of
these events occurs is at most the product of these probabilities
(which is negligible). 

\begin{lemma}\label{lem-negative-correlation}
Define the events $\{A_i \dvtx  i \in I\}$ as above. Then
\[
\mathbb{P}\biggl(\bigcap_{i\in S} A_i \biggr) \leq\prod_{i\in S}\mathbb
{P}(A_i)\qquad \mbox{for any subset $S \subset I$}.
\]
\end{lemma}

\begin{pf}
The proof proceeds by induction on $|S|$. For the empty set this is
trivial, and we will prove
that for every set $S$ and $j \in I \setminus S$
\[
\mathbb{P}\biggl(\bigcap_{i\in S} A_i\cap A_j \biggr) \leq\mathbb{P}
\biggl(\bigcap_{i\in S}
A_i \biggr) \mathbb{P}(A_j).
\]
Define the following independent random variables for every time $t$:
\begin{eqnarray*}
\mathbb{P}(B_t = 1) &=& p_{tj},\qquad \mathbb{P}(B_t = 0) = 1-p_{tj},\\
\mathbb{P}(H_t = i) &=& \frac{p_{ti}}{1-p_{tj}}\qquad \mbox{for each $i
\neq j$}.
\end{eqnarray*}
We may now define $J_t$, the position of the ball at time $t$, as a
function of $B_t$ and~$H_t$,
such that indeed $\mathbb{P}(J_t = i) = p_{ti}$ for all $i$:
\[
J_t = \cases{
j & \quad $B_t = 1$,\vspace*{2pt}\cr
H_t & \quad $B_t = 0$.}
\]
Crucially, the event $A_j$ depends only on the values of $\{B_t\}$, and
is a monotone decreasing in them.
Further notice that the function
\[
f(b_1,\ldots,b_n) \stackrel{\triangle}{=}\mathbb{P}\biggl(\bigcap
_{i\in S} A_i \mid B_1=b_1,\ldots
,B_n=b_n \biggr)
\]
is monotone increasing in the $b_i$'s. Therefore, applying the
FKG-inequality (see, e.g., \cite{AS}, {Chapter 6}, and also \cite
{Grimmett}, {Chapter 2}) on $Y = f(B_1,\ldots,B_n)$ and $\bolds
{1}_{A_j}$ gives
\[
\mathbb{P}\biggl(\bigcap_{i\in S} A_i\cap A_j \biggr)= \mathbb{E}[Y
\bolds{1}_{A_j}] \leq\mathbb{E}
[Y] \mathbb{P}(A_j) = \mathbb{P}\biggl(\bigcap_{i\in S} A_i \biggr)
\mathbb{P}(A_j),
\]
as required.
\end{pf}
Altogether, we obtain that the probability that all of the bins with
indices in $I$
have at most $(\frac{\log n}{\log\log n} - 4)$ balls is
\[
\mathbb{P}\biggl(\bigcap_{i\in I} A_i\biggr) \leq\biggl(1-\frac{ \log^4 n}{
n}\biggr)^{n/\log^2 n} \leq\mathrm{e}^{-\log^2 n}.
\]
This completes the proof of part~(i).
\end{pf*}

\begin{pf*}{Proof of Theorem~\ref{thm-nonadaptive} [\normalfont{Part~(ii)}]} 
It is convenient to describe the proof of this part
for a slightly different model instead of the one considered
in the previous sections. Namely, in the variant model, in every round
each bin among the $n$ bins is chosen randomly and independently
as one of the
options with probability $\alpha$. By Chernoff's bounds, our results
in this model will carry into
the original one, since obtaining $\alpha n$ uniform bins is dominated
by getting each bin independently with probability $(1+\varepsilon
)\alpha
$, and dominates a probability of $(1-\varepsilon)\alpha$ for each bin.

For the simplicity of the notations, we will henceforth consider the
case \mbox{$k=n/2$}, noting that our proofs hold for $k=\alpha n$ with any $0
< \alpha< 1$ fixed.

As noted in the \hyperref[sec1]{Introduction}, the relaxed
model of strategies $Q_t$ such that $\|Q_t\|_\infty\leq k/n$
is stronger than the model where there are $k$ uniform options for bins.
In fact, the results of this part (an optimal maximal load of order
$\sqrt{\log n}$)
do not hold for the relaxed model. For instance, if $Q_t$ assigns probability
$k/n = \frac12$ to $i=t$ and $i=(t+1)$ (with the indices reduced
modulo $n$),
the maximum load will be at most $2$.


However, it is easy to see
that in fact each strategy $Q_t$ is more restricted.
Indeed, the total probability that $Q_t$ can assign to any $r$ bins
does not
exceed $1-2^{-r}$, as for each fixed set $I$ of $r$ bins,
the probability that none of the members of $I$ is an optional
choice for ball number $t$ is $2^{-r}$.

We start with the simple
proof of the upper bound, obtained by the natural algorithm which places
the ball in round $t$ in the first possible bin
(among the $k$ given choices)
that follows bin number $t$ in the cyclic
order of the bins.

\begin{lemma}
\label{l62}
There exists a nonadaptive strategy ensuring that, w.h.p., the maximum
load in the above model is
at most $O(\sqrt{\log n})$.
\end{lemma}

\begin{pf}
Order the bins cyclically $b_1,b_2, \ldots,b_n,b_{n+1}=b_1$.
For each round~$t$, $1 \leq t \leq n$, place the ball number $t$
in the first possible bin $b_i$ that follows $b_t$ in our cyclic order
and is one of the given options for this round. Note, first, that
the probability that the ball in round $t$ is placed in a bin whose
distance from $b_t$ exceeds $2 \log n$, is precisely the
probability that none of the $2 \log n$ bins following $b_t$ is
chosen in round $t$, which is
\[
2^{-2 \log n} < n^{-5/4}.
\]
Therefore, with high probability, this does not happen for any $t$.
In addition, the probability that a fixed bin $b_i$
gets a load of $ \sqrt{\log n}$ from balls placed in the
$2 \log n- 2 \sqrt{\log n}$
rounds $\{i-2 \log n +1, i-2 \log n+2 , \ldots,i-2 \sqrt{\log n}\}$,
does
not exceed
\[
\pmatrix{ 2 \log n -2 \sqrt{\log n}\vspace*{3pt}\cr \sqrt{\log n}}
\biggl(\frac{1}{2^{2 \sqrt{\log n}}}\biggr)^{\sqrt{\log n}}
\leq\frac{1}{n^{2-o(1)}}.
\]
Indeed, for each fixed value of $t \in[i-2 \log n+1,
i-2 \sqrt{\log n} ]$, if the ball placed in
round number $t$ ends in bin number $i$, then none of the
$2 \sqrt{\log n}$ bins preceding $b_i$ is chosen as an optional
bin for ball number $t$, and the probability of this event is
$2^{-2 \sqrt{\log n}}$. There are ${2 \log n -2 \sqrt{\log n}\choose \sqrt{\log n}}$ possibilities to select
$ \sqrt{\log n}$ rounds in the set
$\{i-2 \log n +1, i-2 \log n+2 , \ldots,i-2 \sqrt{\log n}\}$,
and as the choices of options for each round are independent, the
desired estimate follows.

We conclude that with high probability no bin $b_i$ gets any
balls from round $t$ with $t \leq i-2 \log n$, and no bin $b_i$ gets
more than $\sqrt{\log n}$ balls from rounds $t$ with
$i-2 \log n < t \leq i-2\sqrt{\log n}$. As $b_i$ can get at
most $2 \sqrt{\log n}$ balls from all other rounds $t$, (as there
are only $2 \sqrt{\log n}$ such rounds), it follows
that with high probability the maximum load does not exceed $3
\sqrt{\log n}$, completing the proof of the lemma. Note that it is
easy to improve the constant factor $3$ in the estimate proved
here, but we make no attempt to optimize it.
\end{pf}

We proceed with the proof of the lower bound. As in the proof of
part~(i) of the theorem,
let $P=(p_{it})$ be the $n\times n$ matrix of probabilities
corresponding to our nonadaptive strategy, where
$p_{it}$ is the probability that
the ball in round $t$ will be placed in bin number $i$. Recall that
for each
fixed round $t$, the sum of the largest $r$ numbers $p_{it}$ cannot
exceed $1-2^{-r}$.
This fact will be the only property of the distribution
$p_{it}$ used in the proof.

Call an entry
$p_{it}$ of the matrix $P$ \textit{large} if
$p_{it} \geq2^{-\sqrt{\log n}}$, otherwise, $p_{it}$ is \textit{small}. Call a column $t$
of $P$ \textit{concentrated} if it has at least
$\frac{\sqrt{\log n}}{2}$ large elements. We consider two possible
cases.

\begin{itemize}
\item \textit{Case 1}:  There are at least $n/2$ concentrated columns.
\end{itemize}

In this case, there are at least $\frac{n \sqrt{\log n}}{4}$ large
entries in $P$. If there is a row, say row number $i$
of $P$, containing at least, say,
$2^{2 \sqrt{\log n}}$ large entries, then the expected number of
balls in the corresponding bin is $\sum_{t=1}^n p_{it} > 2^{\sqrt
{\log n}} $, and, as the variance of this quantity is
smaller than the expectation, it follows that in this case with high
probability this bin will have a load that exceeds
$\Omega( 2^{\sqrt{\log n}}) > \sqrt{\log n}$. We thus assume that no row
contains more than $2^{2 \sqrt{\log n}}$ large elements.
Therefore, there are at least $\frac{n}{2^{2 \sqrt{\log n}}}
=n^{1-o(1)}$ rows,
each containing at least, say, $\frac{\sqrt{\log n}}{8}$ large
elements. Indeed, we can select such rows one by one. As long as the
number of selected rows does not exceed $\frac{n}{2^{2 \sqrt
{\log}}}$, the total number of large elements in them is\vspace*{-2pt} at most
$n$, and hence the remaining rows still contain at least
$\frac{n \sqrt{\log n}}{4}-n > \frac{n \sqrt{\log n}}{8}$ large
elements, implying that there is still another row containing at
least $\frac{\sqrt{\log n}}{8}$ large elements.

Fix a bin corresponding to a
row with at least $\frac{\sqrt{\log n}}{8}$ large
elements, and fix $\frac{\sqrt{\log n}}{8}$ of them.
The probability that all balls corresponding to these
large elements\vspace*{1pt} will be placed in this bin is at least
\[
\biggl(\frac{1}{2^{\sqrt{\log n}}}\biggr)^{{\sqrt{\log n}}/{8}}
=\frac{1}{n^{1/8}}.
\]
As there are at least $n^{1-o(1)}$ such bins, and the events of
no large load in distinct bins are negatively correlated (see
Lemma~\ref{lem-negative-correlation}), we
conclude that the probability that none of these bins has a load
at least $\frac{\sqrt{\log n}}{8}$ is at most
\[
\biggl(1-\frac{1}{n^{1/8}}\biggr)^{n^{1-o(1)}}=o(1),
\]
showing that in this case
the maximum load is indeed $\Omega(\sqrt{\log n})$ with high
probability.

\begin{itemize}
\item \textit{Case 2}:  There are less than $n/2$ concentrated columns.
\end{itemize}

In this case, the sum of all small entries of the matrix $P$ is at
least
\[
\frac{n}{2 \cdot2^{0.5 \sqrt{\log n}}},
\]
since each of the $n/2$ nonconcentrated columns has less
than $0.5 \sqrt{\log n}$ large elements, and hence in each such
column the sum of all small elements is at least
$2^{-0.5 \sqrt{\log n}}$.

Call a small entry $p=p_{ij}$ of $P$ an entry of \textit{type} $r$
(where $\sqrt{\log n} \leq r \leq2 \log n$), if
$\frac{1}{2^{r+1}} \leq p < \frac{1}{2^r}$. Since the sum of all
entries of $P$ that are smaller than \mbox{$2^{-2 \log n} =1/n^2$} is at
most $1$, there is a value of $r$ in the above range, so that the
sum of all entries of $P$ of type $r$ is at least
\[
\frac{n}{4 \log n \cdot2^{0.5 \sqrt{\log n}}}
> \frac{n}{2^{0.75 \sqrt{\log n}}}.
\]
Put
\[
x\stackrel{\triangle}{=}2^{0.75 \sqrt{\log n}},
\]
and note that there are at least
$\frac{n2^r}{x}$ entries of type $r$ in $P$ (since otherwise their
total sum cannot be at least $n/x$). We now restrict our attention
to these entries.

We can assume that there is no row containing more than
$2^{r+1} \log n$ of these entries. Indeed, otherwise the expected
number of balls in the corresponding bin is at least $\log n$, the
variance is smaller, and hence by Chebyshev with high probability
the load in this bin will exceed
$\Omega(\log n) >\sqrt{\log n}$. We can now apply
again\vspace*{1pt} our greedy procedure and conclude that there are at least
$\frac{n}{4x \log n} =n^{1-o(1)}$ rows, each containing at least
$\frac{2^r}{2x}$ entries of type $r$; indeed, a set of less
than $\frac{n}{4x \log n}$ rows contains a total of at most
\[
\frac{n}{4x \log n} 2^{r+1} \log n =\frac{n 2^r}{2x}
\]
elements
of type $r$, leaving at least $\frac{n 2^r}{2x}$ such elements in
the remaining rows, and hence ensuring the existence of an
additional row with at least $\frac{2^r}{2x}$ such entries.

Fix a bin corresponding to a
row with at least $\frac{2^r}{2x}$ entries of type $r$.
The probability that exactly $t$ balls corresponding to these entries
will be placed in this bin is at least
\[
\pmatrix{ \dfrac{2^r}{2x}\vspace*{3pt}\cr t} \biggl(\frac{1}{2^{r+1}}\biggr)^t
\biggl(1-\frac{1}{2^r}\biggr)^{{2^r}/{(2x)}}
\geq
\biggl(\frac{2^r}{2xt}\biggr)^t \biggl(\frac{1}{2^{r+1}}\biggr)^t
\biggl(1-\frac{1}{2x}\biggr)
> \frac{1}{2}(4xt)^{-t}.
\]
for $t=\sqrt{\log n}$ the last quantity is at least
\[
\tfrac{1}{2} \bigl(4 \cdot2^{0.75
\sqrt{\log n}} \sqrt{\log n}\bigr)^{-\sqrt{\log n}}
=n^{-3/4-o(1)}.
\]
This, the fact that there are $n^{1-o(1)}$ such rows, and the
negative correlation implies that in this case, too, with high
probability there is a bin with load at least
$\Omega(\sqrt{\log n})$. This completes the proof of part~(ii) of the theorem.
\end{pf*}

\section{Concluding remarks and open problems}\label{sec:conclusion}

\begin{itemize}
\item We have established a sharp choice-memory tradeoff for achieving
a constant maximal load in the balls-and-bins experiment,
where there are $n$ balls and $n$ bins, each ball has $k$ uniformly
chosen options for bins, and there are $m$ bits of memory available. Namely:
\begin{enumerate}[1.]
\item[1.] If $k m = \Omega(n)$ for $k = \Omega(\log n)$ and $m = \Omega
(\log n \log\log n)$, then there exists an algorithm that achieves an
$O(1)$ maximal load w.h.p.
\item[2.] If $k m = o(n)$ for $m = \Omega(\log n)$, then any algorithm
w.h.p. creates an unbounded maximal load. For this case we
provide two lower bounds on the load: $\Omega(\log\log
(\frac{n}{km}))$ and
$(1+o(1))\frac{\log(n/m)}{\log\log(n/m)+\log k}$.
\end{enumerate}

\item In particular, if $m = n^{1-\delta}$ for some $\delta> 0$ fixed
and $2 \leq k \leq\operatorname{polylog}(n)$, we obtain
a lower bound of $\Theta(\frac{\log n}{\log\log n})$ on the maximal
load. That is, the typical maximal load in \textit{any} algorithm has
the same order as the typical maximal load in a \textit{random}
allocation of $n$ balls in $n$ bins.

\item Given our methods, it seems plausible and interesting to improve
the above lower bounds to $(1+o(1))\frac{\log({n}/{(km)})}{\log
\log({n}/{(km)})}$, analogous to the load of $(1+o(1))\frac{\log
n}{\log\log n}$
in a completely random allocation.

\item Note that, when $k m = n^{1-\delta}$ for some fixed $\delta> 0$,
even the above conjectured lower bound is still a factor of $\delta$
away from the upper bound given by a random allocation.
It would be interesting to close the gap between these two bounds.
Concretely, suppose that $k m = \sqrt{n}$; can one outperform
the typical maximal load in a random allocation?

\item To prove our main results, we study the problem of achieving a
perfect allocation (one that avoids collisions, i.e., a matching) of
$(1-\delta)n$ balls into $n$ bins. We show that there exist constants
$C > c > 0$ such that:
\begin{enumerate}[1.]
\item If $k m > C n$ for $k = \Omega(\log n)$ and $m = \Omega(\log n
\cdot\log\log n)$, then there exists an algorithm that achieves
a perfect allocation w.h.p.
\item If $k m < c n$ for $m = \Omega(\log n)$, then any algorithm
creates $\Omega(n)$ collisions w.h.p.
\end{enumerate}

\item In light of the above, it would be interesting to show that there
exists a critical $c > 0$ such that, say for $k,m \geq\log^2 n$, the
following holds: If $k m \geq(c+o(1))n$ then there is an algorithm
that achieves a perfect allocation w.h.p., whereas if
$k m \leq(c-o(1))n$ then any algorithm has $\Omega(n)$ collisions
w.h.p.

\item The key to proving the above results is a combination of
martingale analysis and a Bernstein--Kolmogorov type large deviation
inequality. The latter, Proposition~\ref{prop-predictions}, relates a
sum of a sequence of random variables to the sum of its conditional
expectations, and crucially does \textit{not} involve the length of the
sequence. We believe that this inequality may have other applications
in combinatorics and the analysis of algorithms.

\item We also analyzed the case of nonadaptive algorithms, where we
showed that for every $k = O( n \frac{\log\log n}{\log n}
)$, the best
possible maximal load w.h.p. is $\Theta(\frac{\log
n}{\log\log
n})$, that is, the same as in a random allocation. For $k = \alpha n$
with $0 < \alpha< 1$, we proved that the best possible maximal load is
$\Theta(\sqrt{\log n})$. Hence, one can ask what the minimal order of
$k$ is, where an algorithm can outperform the order of the maximal load
in the random allocation.
\end{itemize}

\section*{Acknowledgments}

We thank Yossi Azar and Allan Borodin for helpful discussions, as well
as Yuval
Peres for pointing us to the reference for Theorem~\ref{thm-martingale}.
We also thank Itai Benjamini for proposing the problem of balanced
allocations with limited memory.

%

\printaddresses

\end{document}